\documentclass[12pt, reqon]{amsart}

\usepackage{cite}

\usepackage{color}
 \usepackage{amscd}

\usepackage{amsrefs}

 \usepackage{amssymb, amsmath}

 \input xy
  \xyoption{all}
 
 \input diagxy

\usepackage{epsfig}

    \newtheorem{theorem}{Theorem}[section]
\newtheorem{prop}{Proposition}[section]
\newtheorem{lemma}[theorem]{Lemma}

  \numberwithin{equation}{section}
 \numberwithin{prop}{section}
\newtheorem{exmp}{Example}[section]

\newcommand{\Beq}{\begin{equation}}
\newcommand{\Eeq}{\end{equation}}
\newcommand{\Beqr}{\begin{eqnarray}}
\newcommand{\Eeqr}{\end{eqnarray}}

\newcommand{\rmc}{{\mathbf {Cat}}}

\newcommand{{\wtlg}}{\widetilde\gamma }

\newcommand{{\wtlG}}{\widetilde\Gamma }

\newcommand{{\wtlv}}{\widetilde v }

\newcommand{{\tlg}}{\tilde\gamma }

\newcommand{{\tlG}}{\tilde\Gamma }

 \newcommand{\mbl}{{{\mathbf \Theta}}}

 \newcommand{\mbeta}{{{\mathbf \Psi}}}

\newcommand{\mbg}{\mathcal G}

\newcommand{\mbps}{\mathbf {\Phi} }

\newcommand{\Obj}{{\rm Obj }}
\newcommand{\Mor}{{\rm Mor }}

\begin{document}
\title{On ${\mathbf {Cat}}$-valued sheaves}

\author{Saikat Chatterjee }
\address{School of Mathematics, Indian Institute of Science Education and Research\\
CET Campus\\ Thiruvananthapuram, Kerala-695016\\
India}
\email{saikat.chat01@gmail.com}

\keywords{Presheaves; Sheaves; Union of subcategories; Categorical groups}
\subjclass[2010]{Primary 18F20;  Secondary:18F99 }


\def\xypic{\hbox{\rm\Xy-pic}}

\begin{abstract}
 Let ${\widetilde {\mathcal O}}(\mathbf B)$ be the category of (open) subcategories of a topological groupoid ${\mathbf B}.$ This paper concerns with the ${\mathbf {Cat}}$-valued  sheaves over category ${\widetilde {\mathcal O}}(\mathbf B).$ Since ${\mathbf {Cat}}$ is not a concrete category, traditional definition of presheaf can not deal with the situation. [13] proposes a new framework for the purpose.  Starting from the definition given in [13], we build-up the frame work for ${\mathbf {Cat}}$-valued  sheaves. For that purpose we  introduce a notion of categorical union, such that categorical union of subcategories is a subcategory, which is required for a meaningful definition of a categorical cover of a topological category. The main result is the following. For a fixed category $\mathbf C,$ the categories of local functorial sections  from $\mathbf B$ to $\mathbf C$ define a ${\mathbf {Cat}}$-valued sheaf on ${\widetilde {\mathcal O}}(\mathbf B).$ Replacing $\mathbf C$ with a categorical group $\mathcal G,$ we find a ${\mathbf {CatGrp}}$-valued sheaf on ${\widetilde {\mathcal O}}(\mathbf B).$ 
\end{abstract}



\maketitle

\section{Introduction}\label{s:int}
This paper is the second part of the sequel. First part [13] deals with ${\mathbf {Cat}}$-valued presheaves  and ${\mathbf {Cat}}$-valued sieves. For this paper we will follow-up the definition of  ${\mathbf {Cat}}$-valued presheaves introduced in [13] and develope the notion of ${\mathbf {Cat}}$-valued sheaves. We have reviewed and recalled all the required results (without proofs) of  [13] here as well, so this paper can be read as a more or less   self contained piece. Before we get into the topic of this paper, let us recollect the thoughts and motivation, explained in Section 1 of [13], behind introducing  ${\mathbf {Cat}}$-valued presheaves.

Traditionally a   presheaf  is defined, for a given category $\mathbf C,$ as a contravariant functor [9, 21, 23]
\begin{equation}\label{E:introtradpre}
	R:{\mathbf C}^{\rm op}\longrightarrow {\mathbf {Set}},
\end{equation}
where ${\mathbf {Set}}$ is the (locally small) category of small sets. ${\mathbf C}$ is typically chosen to be the category ${\widetilde {\mathcal O}}(B)$ of open subsets of a topological space $B$, that is, 
\begin{equation}\label{E:introtildeo}
\begin{split}
&\Obj\biggl({\widetilde {\mathcal O}}(B)\biggr) := \{U|\hskip 0.15 cm U\subset B\},\\
&{\rm Hom}(U, V) := \{f:U\to V|\hskip 0.15 cm U, V\subset B\}.
\end{split}
\end{equation}
In \eqref{E:introtradpre} instead of ${\mathbf {Set}},$ one may possibly consider any other \textit{concrete category} such as category of (small) groups ${\mathbf {Grp}}$, category of (small) vector spaces ${\mathbf {Vect}}$ or category of (small) rings ${\mathbf {Ring}},$ and accordingly  respectively gets \textit{presheaf of groups, presheaf of vector spaces} or \textit{presheaf of rings}. We will  often adopt an alternate terminology for  them, such as,   presheaf of groups will be called ${\mathbf {Grp}}$-valued presheaf, and like wise
 \begin{itemize}
 \item{
 presheaf of vector spaces= ${\mathbf {Vect}}$-valued presheaf 
 }
 \item{
 presheaf of rings=${\mathbf {Ring}}$-valued presheaf
 }
\end{itemize} 
and so on.

Now  suppose instead of a topological space $B,$ we are interested in a \textit{topological category} $\mathbf B.$ By a topological category we mean a category  whose both object and morphism spaces are topological spaces. Though this taxonomy of topological category is not universal in literature. We adopt the definition of [11]. In this context natural object of interest  would be the category ${\widetilde {\mathcal O}}(\mathbf B)$ of subcategories of $\mathbf B.$  Natural choice would be to consider the category ${\mathbf {Cat}}$ of small categories as the codomain  of a presheaf in this context, rather than ${\mathbf {Set}}$ (or any other concrete category). However ${\mathbf {Cat}}$ is not a concrete category, and we can not proceed with the existing definition of presheaf given in \eqref{E:introtradpre},  and we need a new framework. In [13] we propose such a framework for ${\mathbf {Cat}}$-\textit{valued presheaves}, and develop  corresponding theory of ${\mathbf {Cat}}$-\textit{valued sieves}. In this paper we establish the corresponding notion of sheaves, namely ${\mathbf {Cat}}$-\textit{valued sheaves}.

We should take note of a crucial issue here. ``\textit{Categorical intersection}" of two subcategories, defined simply as intersection of objet sets and morphism sets, is a subcategory. However, unlike subsets, union of two subcategories, defined naively as union of objet sets and morphism sets, is not a subcategory. But, in order to define a sheaf over a topological category, we need a ``reasonable definition"  of a  cover of a topological category, and in turn that requires a ``reasonable definition" of union of (open) subcategories. So, before progressing to  ${\mathbf {Cat}}$-valued sheaves from ${\mathbf {Cat}}$-valued presheaves, we work out a definition of ``\textit{categorical union}" of subcategories. The categorical unions and categorical intersections satisfy usual inter-relations enjoyed by their set theoretic counterparts. We also define \textit{open subcategories}, and \textit{open categorical cover} of an open subcategory. We construct an example of ``${\mathbf {Cat}}$\textit{-valued sheaf of functorial sections}'' for a fixed category $\mathbf C,$ where an object ${\mathbf U}\in \Obj({\widetilde {\mathcal O}(\mathbf B)})$ is sent to the  category of functors from $\mathbf U$ to $\mathbf C$. Here things get more interesting! Since categorical union of subcategories is not merely union of morphism sets and object sets, it is really a ``non-trivial'' task to establish that categories of local functorial sections indeed define a sheaf (that is, they satisfy appropriate ``locality'' and ``gluing" conditions). In fact at first sight it may seem that this construction is not going to work. However, it is remarkable that despite all the intricacies, it turns out that local functorial sections (and natural transformations between them) do define a ${\mathbf {Cat}}$-valued sheaf.
We also consider the \textit{sheaves of categorical groups} or ${\mathbf {CatGrp}}$-\textit{valued sheaves}, where ${\mathbf {CatGrp}}$ is the category of {categorical groups}, treated as a full subcategory (identified via an obvious full faithful forgetful functor ${\mathbf {CatGrp}}\longrightarrow{\mathbf {Cat}}$) of ${\mathbf {Cat}}.$ 

Let $\mathcal C$ be a category of a collection small categories. We work with a ${\mathbf {Cat}}$-\textit{valued presheaf} over $\mathcal C$ given by a contravariant functor:
\begin{equation}\label{E:introtradprecat}
{\mathcal R}: {\mathcal C}^{\rm op}\longrightarrow {\mathbf {Cat}}.
\end{equation}
The back ground motivation for our construction stems from the study of categorical geometry [8,14,15,18,34] on a \textit{path space groupoid} of a given smooth manifold $M$; that is, a category  ${\mathbb P}M$, whose object space is the manifold $M$ and morphisms are certain equivalence classes of smooth paths. Usual compact-open topology defines a topology on $\Mor({\mathbb P}M).$ The path space groupoid over a smooth manifold naturally occurs providing the back-ground geometry  in higher gauge theories [7,22,28,30,33] and non abelian gerbe theories [1--6, 10, 31, 32]. Recently the ``locally defined subcategories" of ${\mathbb P}M$ also appeared in the context of local trivializations of categorical principal bundles [15, 19]. In the context of this paper, it would be of particular interest to consider the case when ${\mathcal C}={\widetilde {\mathcal O}}({\mathbb P}M)$ . However considering the length of this paper, we do not pursue  a detailed and rigorous treatment of ${\mathbb P}M.$ For that reason, though the framework  developed here is perhaps most suitable for ${\mathcal C}={\widetilde {\mathcal O}}({\mathbb P}M),$ we tread  on a more abstract approach and try to minimize the reference to  ${\mathbb P}M.$ Occasionally when we have to recall  ${\mathbb P}M,$ we would give a heuristic description without going into the technicalities. For a detailed exposition on the topic, we refer to [5, 16, 29].
 \subsection*{Notations and conventions}
 We borrow our notations from [13]. Let $\mathbf C$ and $\mathbf D$ be a given pair of categories.
\begin{equation}\label{E:notfunc}
	{\rm Fun}(\mathbf C, \mathbf D)
\end{equation}
will denote the set of all functors from $\mathbf C$ to $\mathbf D.$ The set of all natural transformations between functors from $\mathbf C$ to $\mathbf D$ will be denoted as
\begin{equation}\label{E:notnat}
	{\mathcal N}(\mathbf C, \mathbf D).
\end{equation}
If ${\mathbf \theta}_1, {\mathbf \theta}_2:\mathbf C\to \mathbf D$ are a pair of functors, then 
\begin{equation}\label{E:notnatspfunc}
	{\rm Nat}(\mathbf \theta_1, \theta_2)
\end{equation}
is the set of all natural transformations between $\theta_1$ and $\theta_2.$  We will often denote a natural transformation $\Phi$ from a functor $\theta_1$ to another functor $\theta_2$ as
\begin{equation}
\Phi:\theta_1\Longrightarrow \theta_2.
\end{equation}
The category of functors will be denoted as
\begin{equation}\label{E:notfunccat}
	{\mathcal F}(\mathbf C, \mathbf D);
\end{equation}
that is 
\begin{equation}\label{E:notfunccatobjmor}
	\begin{split}
		&\Obj\biggl({\mathcal F}(\mathbf C, \mathbf D)\biggr)={\rm Fun}(\mathbf C, \mathbf D)\\
		&\Mor\biggl({\mathcal F}(\mathbf C, \mathbf D)\biggr)={\mathcal N}(\mathbf C, \mathbf D).
\end{split}
\end{equation}
Given a morphism $f$ in some category, $s(f), t(f)$ will respectively denote the source of $f$ and target of $f;$ that is,
$${s(f)}\xrightarrow{f}{t(f)}.$$

$\mathbf {\emptyset}$ will denote the \textit{empty category}; i.e. a category whose object and morphism sets are null sets.

  \subsection*{Overview of the sections} 
  We now give a brief overview of each section of this paper.
  \begin{itemize}
 \item{In Section~\ref{S:mbsieve} we  review some of the constructions and results from [13]. We work with   a category  of a collection of  small categories, $\mathcal C$. In particular, in this section we recall the definition of  a ${\mathbf {Cat}}$-valued presheaf over ${\mathcal C}$ introduced in [13]. We state a version of Yoneda embedding in Proposition~ \ref{Th:yonembed}.   
  
   }
   \item{ Section~\ref{s:cpb} prepares the stage for the next section.  We first define the categorical union of subcategories, and Proposition~\ref{pr:boolean}   ensures that the categorical union and intersection satisfy Boolean relations (Proof carried out in the Appendix).  We take $\mathbf B$ to be a groupoid. We call a subcategory open if  both object and morphism sets are open, and the subcategory is a groupoid.   This particular definition of open subcategories is essential to have a well defined sheaf of functorial sections. We introduce the notion of an open categorical cover and provide some examples for the same.   }

 \item{ Section~\ref{S:catsheaves} is the central theme of this paper. The ultimate goal in this section is to construct the (${\mathbf {Cat}}$-valued) sheaf of functorial sections. 
 We define a ${\mathbf {Cat}}$-valued sheaf   over ${\widetilde {\mathcal O}}(\mathbf B)$ to be a ${\mathbf {Cat}}$-valued presheaf satisfying certain locality and gluing conditions (described in \eqref{E:objsheafcov}--\eqref{E:ovmor}). 
 
 Next we fix a category $\mathbf C.$ Our goal is to construct the (${\mathbf {Cat}}$-valued) sheaf of functorial sections   to category $\mathbf C.$     Finally in Theorem~\ref{Th:funcsheaf} we prove that functorial sections  with respect to category $\mathbf C$ over $\mathbf B$ define a ${\mathbf {Cat}}$-valued sheaf. 
 }
 \item{In section~\ref{S:sheavcaygrop} we first give a brief review of categorical groups, and show (Proposition~\ref{P:funccat}) that ${\mathcal G}^{\mathbf U}:={\mathcal F}(\mathbf U, \mathcal G)$ form a categorical group, where $\mathcal G$ is a fixed categorical group and $\mathbf U$ any category. We define (category of categorical groups) ${\mathbf {CatGrp}}$-valued presheaves. In Proposition~\ref{Pr:funcsheafgrp} we show that if we replace $\mathbf C$ of Section~\ref{s:cpb}   by a categorical group ${\mathcal G},$ then we  get a ${\mathbf {CatGrp}}$-valued sheaf on ${\widetilde {\mathcal O}}(\mathbf B).$

 }
  \item{ The Appendix contains the proof of Proposition~\ref{pr:boolean}.  
  }
  \item{We end this paper with some concluding remarks to highlight future scopes of our construction and results  in this paper.}
  \end{itemize}         
\subsection*{Summary of the paper}  Main objective of this paper is to propose a framework for $\mathbf {Cat}$-valued sheaves and construct an example of  ``${\mathbf {Cat}}$-valued {sheaf of (local) functorial sections}'' defined on  a topological groupoid $\mathbf B$ with respect to  a fixed category $\mathbf C$. Our starting point is the  $\mathbf {Cat}$-valued presheaf on ${\mathcal C}$ introduced in [13] and recalled in 
Section~\ref{S:mbsieve}. 

A $\mathbf {Cat}$-valued presheaf on ${\mathcal C}$ is a contravariant functor
$${\mathcal C}^{\rm op}\longrightarrow {\mathbf {Cat}}.$$

With the motivation to define $\mathbf {Cat}$-valued sheaves on ${\widetilde {\mathcal O}}(\mathbf B),$ we introduce a new notion of categorical union of subcategories, which ensures that union of subcategories is a subcategory. We show that categorical unions thus defined and intersections of subcategories are consistent with the standard set theoretic relation. We turn to define  open subcategories and categorical cover  of an open subcategory of a topological groupoid $\mathbf B$ and illustrate the definitions with several examples. Then we define a   $\mathbf {Cat}$-valued sheaf on ${\widetilde {\mathcal O}}(\mathbf B)$ to be a $\mathbf {Cat}$-valued presheaf satisfying certain ``gluing'' and ``locality'' conditions.

We fix a category $\mathbf C$ and consider the functor categories ${\mathbf C}^{\mathbf U}:={\mathcal F}(\mathbf U, \mathbf C)$ defined by open subcategories $\mathbf U$ of $\mathbf B.$ We show that the prescription
$${\mathbf U}\mapsto {\mathbf C}^{\mathbf U}$$
then produces a $\mathbf {Cat}$-valued sheaf on ${\widetilde {\mathcal O}}(\mathbf B).$ We call them ${\mathbf {Cat}}$-valued sheaf of \textit{functorial sections} on ${\widetilde {\mathcal O}}(\mathbf B).$

As a special case of above construction, taking $\mathbf C$ to be a categorical group $\mathcal G,$ we obtain a $\mathbf {CatGrp}$-valued sheaf on ${\widetilde {\mathcal O}}(\mathbf B).$

\section{$\mathbf {Cat}$-valued presheaves}\label{S:mbsieve}

For the ease of reference, in this section we briefly review and collect some of the relevant material from{}[13]. Some parts would be verbatim copy from {}[13]. We keep the details to a bare minimum, and skip the proofs. We will provide the reference to corresponding Proposition/Theorem numbers in{}[13].

 Let $\mathcal C$ be a category of a collection of  (small)categories; that is objects are a set of (small) categories and morphisms are functors between them. Later we will mostly be dealing with the category ${\widetilde {\mathcal O}}(\mathbf B),$ where $\mathbf B$ is a given category and, 
\begin{equation}\label{E:genob}
	\begin{split}
		&{\rm Obj}\biggl({\widetilde {\mathcal O}}({\mathbf B})\biggr):=\{{\mathbf U}| {\mathbf U}\subset {\mathbf B}\}=\hbox{set of all subcategories of}\hskip 0.2 cm {\mathbf B},\\
	 &{\rm Hom}({\mathbf U}, {\mathbf V})=\{{\mbl}:{\mathbf U}\longrightarrow {\mathbf V}|{\mathbf U}, {\mathbf V}\subset {\mathbf B}\}.	
\end{split}
\end{equation}

We will also work with the category ${{\mathcal O}}(\mathbf B),$ whose objects are same as those of ${\widetilde {\mathcal O}}(\mathbf B);$ but, only morphism  between any two subcategories (objects) is the  inclusion functor, if one is  subcategory  of the other. Otherwise no morphism exists:

\begin{equation}\label{E:ob}
	\begin{split}
		&{\rm Obj}\biggl({\mathcal O}({\mathbf B})\biggr):=\{{\mathbf U}| {\mathbf U}\subset {\mathbf B}\}=\hbox{set of all subcategories of}\hskip 0.2 cm {\mathbf B},\\
	 &{\rm Hom}({\mathbf U}, {\mathbf V})=\{{\mathbf i}:{\mathbf U}\hookrightarrow {\mathbf V}|{\mathbf U}\subset {\mathbf V}\subset {\mathbf B}\},\\ 
		&\hskip 02.4 cm \hbox{where}\hskip 0.2cm {\mathbf i}\hskip 0.2 cm \hbox{is the inclusion functor, and}\\
	 &{\rm Hom}({\mathbf U}, {\mathbf V})=\emptyset, \hbox {if}\hskip 0.2 cm {\mathbf U}\not\subset {\mathbf V}.
\end{split}
\end{equation}

Let $\rmc$ be the category of all (small) categories. Analogous to a contravariant  $\rm Hom$-functor in set theoretic frame work, we have a contravariant functor ${\mathcal F}_{\mathbf U}:{\mathcal C}^{\rm op}\to \rmc,$ corresponding to each $\mathbf U\in {\rm Obj}(\mathcal C),$ as follows.
\begin{eqnarray}
	&{\mathcal F}_{\mathbf U}:&{\rm Obj}(\mathcal C)\to {\rm Obj}(\rmc)\nonumber\\
	       && \mathbf V\mapsto {\mathcal F}(\mathbf V, \mathbf U)\label{E:objfun}\\
	&{\mathcal F}_{\mathbf U}:&{\rm Mor}(\mathcal C)\to {\rm Mor}(\rmc)\nonumber\\
 &&\biggl({\mathbf V}\xrightarrow{\mbl} {\mathbf W}\biggr)\mapsto \biggl({\mathcal F}(\mathbf W, \mathbf U)\xrightarrow{{\mathcal F}_{\mathbf U}(\mbl)}{\mathcal F}(\mathbf V, \mathbf U)\biggr),\label{E:morfun}
\end{eqnarray}
where $\mathbf \Theta$ is a functor from the category $\mathbf V$ to $\mathbf W,$ and \eqref{E:morfun} is defined by following two equations. 
\begin{equation}\label{D:objfunc}
	\begin{split}
		{\mathcal F}_{\mathbf U}(\mbl):\Obj\biggl({\mathcal F}(\mathbf W, \mathbf U)\biggr)&\longrightarrow \Obj\biggl({\mathcal F}(\mathbf V, \mathbf U)\biggr),\\
		 \mbeta&\mapsto \mbeta\mbl,\\
	& \xymatrixcolsep{5pc}
\xymatrix{
	& {\mathbf V} \ar[d]^{\mbeta\mbl} \ar@{-->}[dl]_-{\mbl} \\
	{\mathbf W} \ar[r]^{\mbeta} &{\mathbf U}
}
\end{split}
\end{equation}
and 
\begin{equation}\label{E:defnat}
	\begin{split}
	&{\mathcal F}_{\mathbf U}(\mbl):\Mor\biggl({\mathcal F}(\mathbf W, \mathbf U)\biggr)\longrightarrow \Mor\biggl({\mathcal F}(\mathbf V, \mathbf U)\biggr),\\
&{\mathcal F}_{\mathbf U}(\mbl):{\mathcal N}(\mathbf W, \mathbf U)\longrightarrow {\mathcal N}(\mathbf V, \mathbf U),\\
&\bigl({\mathcal F}_{\mathbf U}(\mbl)\bigr)(\mathcal S):={\mathcal S}\mbl\in {\mathcal N}(\mathbf V, \mathbf U).
	\end{split}
\end{equation}

Following proposition confirms that indeed ${\mathcal F}_{\mathbf U}:{\mathcal C}^{\rm op}\to \rmc$  is a functor.

\begin{prop}\label{Pr:Hompresheaf}
	{\rm [Proposition 2.1, [13]]} Let $\mathcal C$ be a category of a collection of  (small)categories and $\rmc$ be the category of all (small) categories. Then, for each ${\mathbf U}\in \Obj(\mathcal C),$ we have a contravariant functor ${\mathcal F}_{\mathbf U}:{\mathcal C}^{\rm op}\to \rmc$ defined as in \eqref{E:objfun},  \eqref{D:objfunc} and \eqref{E:defnat}.
\end{prop}
Moreover the functors ${\mathcal F}_{\mathbf U}$ are consistent with the Yoneda lemma.		
\begin{prop}\label{Pr:yoneda}
{\rm [Proposition 2.2, {}[13]]}
There exists an isomorphism between  ${\rm Fun}(\mathbf U, \mathbf V)$ and ${\rm Nat}({\mathcal F}_{\mathbf U}, {\mathcal F}_{\mathbf V}):$
\begin{equation}\label{E:yonlemma}
 {\rm Fun}(\mathbf U, \mathbf V)\cong{\rm Nat}({\mathcal F}_{\mathbf U}, {\mathcal F}_{\mathbf V}).
\end{equation}
\end{prop}


\subsection{Presheaves of categories}
 Two prominent directions of enquiry, which immediately emerge out of the definition of presheaves, are \textit{sheaves} and \textit{sieves}. In {}[13] we have introduced the notion of  $\mathbf {Cat}$-valued presheaf to study the $\mathbf {Cat}$-valued sieves. Here we will recall the definition of  $\mathbf {Cat}$-valued presheaf given in {}[13]. In this paper our focus will be $\mathbf {Cat}$-valued sheaves.

As before, let $\mathcal C$ be a category of a collection of (small) categories, and $\rmc$ be the category of all small categories.  We work with following definition of  \textit{presheaves of categories} over $\mathcal C$. 
A   \textit{presheaf of categories} (or, a $\mathbf {Cat}$-\textit{valued presheaf}), over the category $\mathcal C$,  is a functor 
\begin{equation}\label{E:prshfcat}
	{\mathcal R}:{\mathcal C}^{\rm op}\to \rmc.
\end{equation}	
An immediate consequence of  the definition above and Proposition~\ref{Pr:Hompresheaf} is the following. 
\begin{lemma}\label{cor:funcprsf}
{\rm [Corollary 3.1, {}[13]]} For each $\mathbf U\in \Obj(\mathcal C),$ the functor $${\mathcal F}_{\mathbf U}:{\mathcal C}^{\rm op}\to\rmc,$$ in Proposition~\ref{Pr:Hompresheaf} is a presheaf of categories, over the category $\mathcal C.$
\end{lemma}

Let ${{\mathbf {Prsh}}}(\mathcal C, {\mathbf {Cat}}):={\mathcal F}({\mathcal C}^{\rm op}, {\mathbf {Cat}})$ denote the category of $\mathbf {Cat}$-valued presheaves, over the category ${\mathcal C};$ that is,
\begin{equation}\label{E:catprsheaf}
	\begin{split}
		&\Obj\biggl({{\mathbf {Prsh}}}(\mathcal C, {\mathbf {Cat}})\biggr)={\rm Fun}({\mathcal C}^{\rm op}, {\mathbf {Cat}}),\\
  &\Mor\biggl({{\mathbf {Prsh}}}(\mathcal C, {\mathbf {Cat}})\biggr)={\mathcal N}({\mathcal C}^{\rm op}, {\mathbf {Cat}}).
\end{split}
\end{equation}

Using Lemma~\ref{cor:funcprsf} and Proposition~\ref{Pr:yoneda} we prove a version of Yoneda embedding  of  $\mathcal C$ in  ${{\mathbf {Prsh}}}(\mathcal C, {\mathbf {Cat}}).$ 
\begin{prop}\label{Th:yonembed}
	{\rm [Theorem 3.2, {}[13]]}Let $\mathcal C$ be a category of a collection of (small) categories, and $\rmc$ be the category of all small categories. Let ${{\mathbf {Prsh}}}(\mathcal C, {\mathbf {Cat}}):={\mathcal F}({\mathcal C}^{\rm op}, {\mathbf {Cat}})$ be the  category of $\mathbf {Cat}$-valued presheaves over the category ${\mathcal C}.$ Then there exists a full and faithful functor 
	$$\mathcal C\longrightarrow{{\mathbf {Prsh}}}(\mathcal C, {\mathbf {Cat}}).$$
 In other words,  $\mathcal C$ can be identified as a full subcategory of
	${{\mathbf {Prsh}}}(\mathcal C, {\mathbf {Cat}}).$
\end{prop}

Instead of working with the entire category $\mathbf {Cat},$ one can consider a  presheaf of subcategories of $\mathbf {Cat}$. For example, one may define a presheaf of \textit{categorical groups}, over $\mathcal C,$ to be a contravariant functor from $\mathcal C$ to $\mathbf {CatGrp}:$
$${\mathcal R}: {\mathcal C}^{\rm op}\to \mathbf {CatGrp},$$
where $\mathbf {CatGrp}$ is the category of  categorical groups. We will denote \textit{category of presheaves of categorical groups} by $${{\mathbf {Prsh}}}(\mathcal C, {\mathbf {CatGrp}}).$$ In section~\ref{S:sheavcaygrop}  we will construct such an example of presheaf of categorical groups. 
\section{Categorical cover}\label{s:cpb}
In the next section we are going to introduce the notion of $\mathbf {Cat}$\textit{-valued sheaves}. For that purpose, we need a ``reasonable definition" of a cover of a categorical space by subcategories. 
This section will be devoted to establish a notion of union of subcategories (namely ``\textit{categorical union}") of a given category, and \textit{categorical cover} of a ``topological category''.
\subsection{Categorical cover}
Suppose $\mathbf B$ is a category, such that both object and morphism spaces are topological spaces. 
We will call such categories \textit{topological categories} [11]. If the category is also a groupoid, the we call it a \textit{topological groupoid}. For example, a Lie groupoid is a topological groupoid. We will soon see another example of a topological groupoid, which would be more relevant for our purpose, the ``path space groupoid'' of a topological space.  It is easy to see that \textit{intersection of two  subcategories} 
$\mathbf U, \mathbf V,$ defined as:
\begin{equation}\label{E:objmorint}
\begin{split}
\Obj(\mathbf U\cap \mathbf V):=\Obj(\mathbf U)\cap \Obj(\mathbf V),\\
\Mor(\mathbf U\cap \mathbf V):=\Mor(\mathbf U)\cap \Mor(\mathbf V),
\end{split}
\end{equation}
is also a subcategory. However, if we define the union of two subcategories in similar fashion,
\begin{equation}\label{E:objmoruninot}
\begin{split}
\Obj(\mathbf U\cup \mathbf V):=\Obj(\mathbf U)\cup \Obj(\mathbf V),\\
\Mor(\mathbf U\cup \mathbf V):=\Mor(\mathbf U)\cup \Mor(\mathbf V),
\end{split}
\end{equation}
then this union fails to be a subcategory. One can see that  the obstruction, for above union of subcategories to be a subcategory, is precisely the $\Obj(\mathbf U)\cap \Obj(\mathbf V),$ because of the following reason.
Suppose $a\xrightarrow{f_1}b\in \Mor(\mathbf U)$ and $b\xrightarrow{f_2}c\in \Mor(\mathbf V),$ where $b\in \Obj(\mathbf U)\cap \Obj(\mathbf V).$ Then $f_2, f_1\in \Mor(\mathbf U)\cup \Mor(\mathbf V),$ but
$f_2\circ f_1\notin \Mor(\mathbf U)\cup \Mor(\mathbf V).$ Note that if $\mathbf U\cap \mathbf V$ is empty, then we have a well define union of subcategories defied as in \eqref{E:objmoruninot}. 


\begin{figure}[h]
\begin{center}
\includegraphics[height=3.0in,width=5.5in]{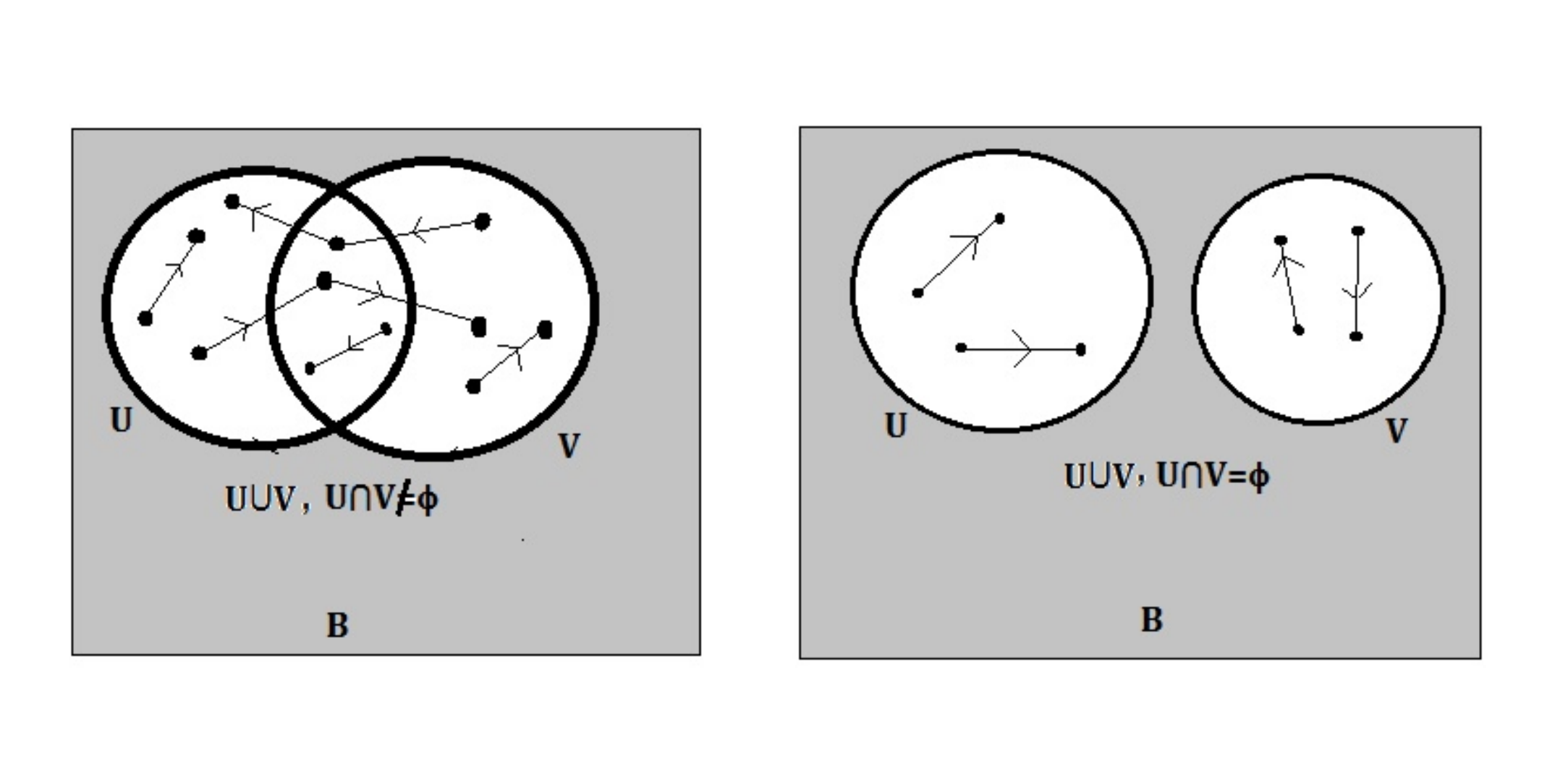}
\caption{A schematic diagram for categorical union of a pair of subcategories}
\end{center}
\end{figure}\label{f:schema}
As a remedy to this problem we define \textit{categorical union} of two subcategories as follows. 
\begin{equation}\label{E:objmoruni}
\begin{split}
&\Obj(\mathbf U\cup \mathbf V):=\Obj(\mathbf U)\cup \Obj(\mathbf V),\\
&\Mor(\mathbf U\cup \mathbf V):={\mathbf {Gen}}\biggl(\Mor(\mathbf U)\cup \Mor(\mathbf V)\biggr),
\end{split}
\end{equation}
where 
\begin{equation}\label{E:spanmor}
\begin{split}
&{\mathbf {Gen}}\biggl(\Mor(\mathbf U)\cup \Mor(\mathbf V)\biggr)=\\
&\bigg\{ f_2\circ f_1|f_2, f_1\in  \Mor(\mathbf U)\cup \Mor(\mathbf V), {s}(f_2)={t}(f_1) \bigg\};
\end{split}
\end{equation}
that is, $\Mor(\mathbf U\cup \mathbf V)$ is the subset of $\Mor(\mathbf B),$ generated by $\Mor(\mathbf U)$ and  $\Mor(\mathbf V).$ Note that $\Mor(\mathbf U)\cup \Mor(\mathbf V)\subset {\mathbf {Gen}}\biggl(\Mor(\mathbf U)\cup \Mor(\mathbf V)\biggr),$ and \eqref {E:spanmor} can be rewritten as
\begin{equation}\label{E:spanmorrewri}
\begin{split}
&{\mathbf {Gen}}\biggl(\Mor(\mathbf U)\cup \Mor(\mathbf V)\biggr)=\\
&\Mor(\mathbf U)\bigcup \Mor({\mathbf V})\bigcup\\
&\Big\{ f_2\circ f_1|f_2 \in  \Mor(\mathbf U), f_1\in\Mor(\mathbf V), {s}(f_2)={t}(f_1) \in \Obj(\mathbf U)\cap \Obj({\mathbf V})\Big\} \bigcup\\ 
&\Big\{ f_2\circ f_1|f_2 \in  \Mor(\mathbf V), f_1\in\Mor(\mathbf U), {s}(f_2)={t}(f_1) \in \Obj(\mathbf U)\cap \Obj({\mathbf V})\Big\}.\\
\end{split}
\end{equation}
For non-intersecting $\mathbf U$ and $\mathbf V,$ we have  $${\mathbf {Gen}}\biggl(\Mor(\mathbf U)\cup \Mor(\mathbf V)\biggr)=\Mor(\mathbf U)\cup \Mor({\mathbf V}),\qquad \hbox{if} \hskip 0.2 cm \mathbf U\cap \mathbf V=\emptyset.$$
It is obvious that the  categorical union of $\mathbf U, \mathbf V,$ defined in \eqref{E:objmoruni}, is a subcategory of $\mathbf B.$ Henceforth by union of subcategories we will always understand a categorical union  of subcategories.

It would be often convenient to visualize  ${\mathbf B} $ as the ``\textit{path space category (groupoid)}" on a (smooth) topological space $B.$  We will not discuss the path space category here in detail. [5, 12, 16] can be consulted for the same. Intuitively, the path space category, over a topological space $B,$ is the category whose objects are points of $B,$ and morphisms are  paths  on $B$ (modulo ``certain equivalence relations"), and composition is given by usual concatenation of  paths. The equivalence relation is typically taken to be ``\textit{thin homotopy equivalence}"  [28] or ``\textit{back-track equivalence}" [16]. The later is slightly weaker condition compare to the former (see Section 3 of [16] for a comparison between the two). Since, every ``path"  can be reversed, this category in fact is a groupoid. The morphism space of a path space groupoid is equipped with usual compact-open topology. We will denote the path space category(groupoid) on $B$ as ${\mathbb P}B,$
\begin{equation}\label{E:pathcat}
\begin{split}
&{\Obj}({\mathbb P}B)=B,\\
&{\Mor}({\mathbb P}B)=\{\hbox{paths on $B$}\}/{\hbox{``some equivalence relations"}}.
\end{split}  
\end{equation}
 Then, we may take $\mathbf U:={\mathbb P}U$ and $\mathbf V:={\mathbb P}V$  to be path space categories on  open subsets $U\subset B$ and $V\subset B$ respectively. Obviously $\mathbf U, \mathbf V$ are subcategories of ${\mathbb P}B := \mathbf B.$ So, in this case  ${\mathbf {Gen}}\biggl(\Mor(\mathbf U)\cup \Mor(\mathbf V)\biggr)$ is basically the set of all paths (modulo ``certain equivalence relations'') lying on $U\cup V\subset B.$ In other words, $\mathbf U\cup \mathbf V$ is the  path space category ${\mathbb P}(U\cup V)$ on $U\cup V\subset B.$ On a related note we observe that even without the smoothness condition on $B$ we can construct the path space groupoid ${\mathbb P}B$ by taking the morphism space to be the space of homotopy equivalent paths on $B$, which is of course stronger equivalence  condition than the previous two. 

The intersection and (categorical) union of subcategories respectively defined in \eqref{E:objmorint}, \eqref{E:objmoruni} satisfy usual set theoretic Boolean relations of unions and intersections. Here we list some of them, and  have provided a proof for the following proposition in the Appendix.
\begin{prop}\label{pr:boolean}
Let $\mathbf U, \mathbf V, \mathbf W$ be subcategories of a category $\mathbf B.$ Suppose intersection and (categorical) union  of subcategories are respectively defined as in \eqref{E:objmorint}, \eqref{E:objmoruni}. Then following relations hold:
\begin{itemize}
\item[]{{\rm (i)} $\mathbf U \cup (\mathbf V \cup\mathbf W) = (\mathbf U \cup \mathbf V )\cup\mathbf W$ and $\mathbf U \cap (\mathbf V \cap\mathbf W) = (\mathbf U \cap \mathbf V )\cap\mathbf W,$}
\item[]{{\rm (ii)}$\mathbf U \cup (\mathbf V \cap\mathbf W) = (\mathbf U \cup \mathbf V )\cap({\mathbf U}\cup{\mathbf W}).$ }
\item[]{{\rm (iii)} $\mathbf U \cap (\mathbf V \cup\mathbf W) = (\mathbf U \cap \mathbf V )\cup({\mathbf U}\cap{\mathbf W}),$         }
\end{itemize}

\end{prop}
\begin{proof} See the Appendix.
\end{proof}

 Let $\mathbf U$ be a sub category  of a topological groupoid $\mathbf B,$ such that $\Obj{(\mathbf U)}$ and $\Mor{(\mathbf U)}$ are open (respectively in $\Obj(\mathbf B)$ and $\Mor(\mathbf B)$).  Then we call $\mathbf U$ to be an \textit{open subcategory} of $\mathbf B,$ if $\mathbf U$ is also a groupoid. 
  The following obvious result would be useful later.
\begin{lemma}\label{L:property}
Suppose $\mathbf U$ is an open  subcategory of a topological groupoid $\mathbf B,$ and for $a\xrightarrow {f_1}b, b\xrightarrow {f_2}c, a\xrightarrow {f'_1}b', b'\xrightarrow {f'_2}c\in \Mor(\mathbf U),$  we have 
\begin{equation}\label{E:samecompo}
f_2\circ f_1=f'_2\circ f'_1.
\end{equation}
Then there exists an isomorphism $b\xrightarrow{g_0}b'\in \Mor(\mathbf U),$ such that 
\begin{equation}\label{E:compocon}
\begin{split}
&g_0\circ f_1=f'_1\\
&f_2\circ g_0^{-1}=f'_2.\hskip 5 cm \fbox{}
\end{split}
\end{equation}
\end{lemma}
\begin{proof}
Directly follows from the fact that $\mathbf U$ is a groupoid. 
\end{proof}

\begin{figure}[h]
\begin{center}
\includegraphics[height=1.67in,width=3.5in]{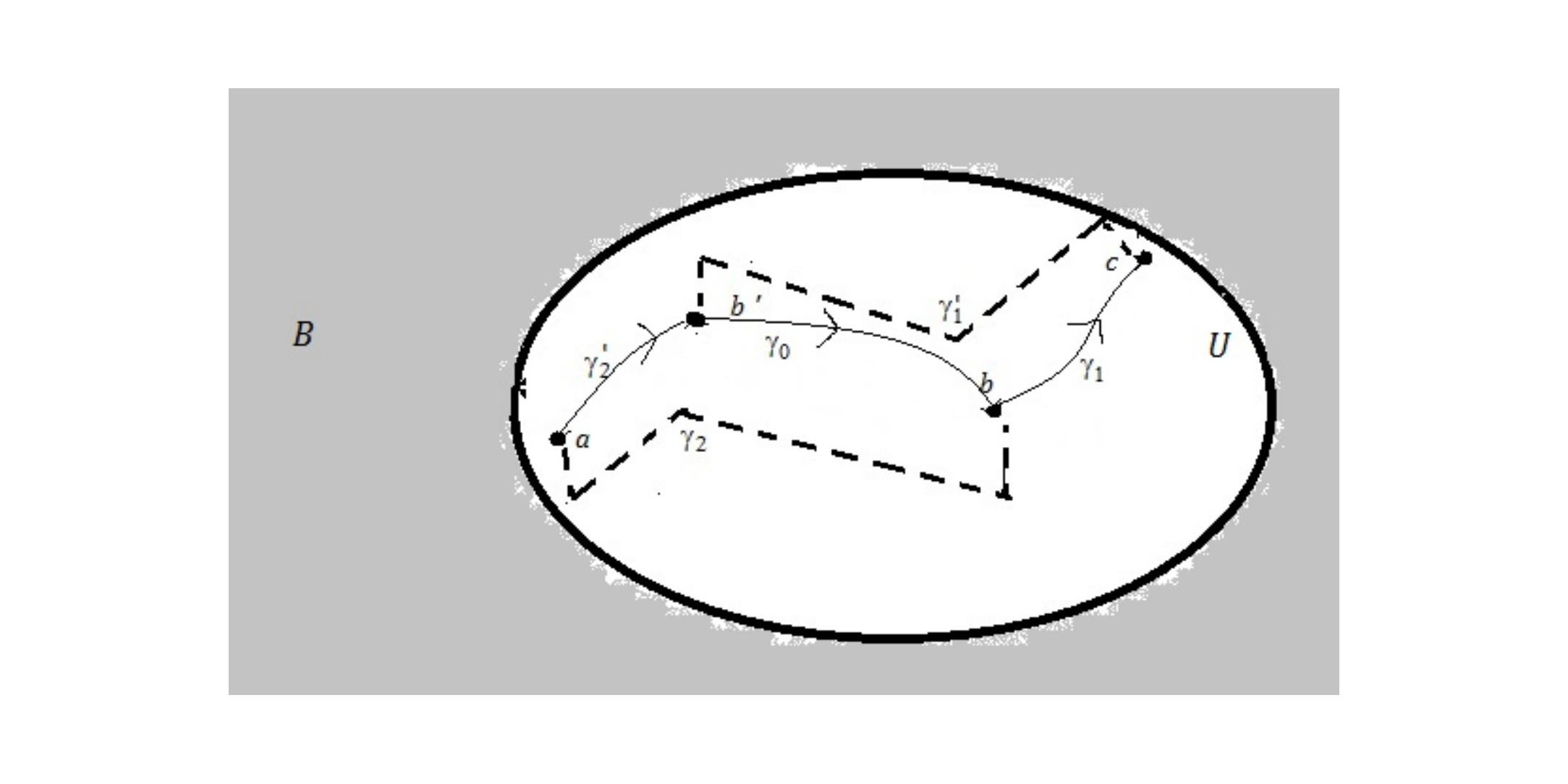}
\caption{}
\end{center}
\end{figure}\label{f:backtrack2}
In the context of the path space groupoid ${\mathbb P}B$ over a smooth topological space $B,$ there is an elegant approach, namely ``\textit{back-track erasing}" [16, 25], to deal with such a situation and the condition in \eqref{E:samecompo}--\eqref{E:compocon} is ensured by imposing certain equivalence relation on the space of paths. Roughly the idea is as follows. If  $U$ is an open subset of $B,$ and  ${\mathbf U}:={\mathbb P}U$ is the path space category over $U$ , that is
\begin{equation}\label{E:pathcatsubset}
\begin{split}
&{\Obj}({\mathbb P}U)=U,\\
&{\Mor}({\mathbb P}U)=\{\hbox{paths on $U$}\}/{\hbox{``some equivalence relations"}}.
\end{split}  
\end{equation}

  Suppose two pairs of paths (modulo equivalence), $a\xrightarrow {\gamma_1}b, b\xrightarrow {\gamma_2}c \in \Mor({\mathbb P}U)$, and $a\xrightarrow {\gamma'_1}b', b'\xrightarrow {\gamma'_2}c\in \Mor({\mathbb P}U)$ such that
$$\gamma_2\circ \gamma_1=\gamma'_2\circ \gamma'_1,$$
then we have the path segment $b\xrightarrow{\gamma_0} b'$ between $b$ and $b',$ satisfying
\begin{equation}\label{E:compoconpath}
\begin{split}
&\gamma_0\circ \gamma_1=\gamma'_1\\
&\gamma_2\circ \gamma_0^{-1}=\gamma'_2.
\end{split}
\end{equation}
Figure 2 illustrates \eqref {E:compoconpath}.

For  future reference we make  following observation.
\begin{lemma}\label{L:uniintprop}
Let $\mathbf U, \mathbf V$ be open subcategories of a topological groupoid $\mathbf B$ . Then 
${\mathbf U}\cup {\mathbf V}$ and ${\mathbf U}\cap {\mathbf V}$ are also open subcategories of $\mathbf B$.
\begin{proof}
Straightforward verification.
\end{proof}
\end{lemma}

Let $\mathbf U$ be an open subcategory of a topological groupoid  $\mathbf B.$ An indexed family of  subcategories $\{{\mathbf U}_{\alpha}\}_{\alpha\in I}$ is an \textit{open categorical cover} of $\mathbf U,$ if each ${\mathbf U}_{\alpha}$ is an open subcategory of  $\mathbf U,$ and
\begin{equation}\label{E:covobjmor}
\mathbf U=\bigcup_{\alpha}{\mathbf U}_{\alpha},
\end{equation}
where the union of subcategories is taken according to the definition in \eqref{E:objmoruni}.

Following examples illustrate  the idea behind an open categorical cover of a topological category.
\begin{exmp}\label{Ex:trivcov}
\rm{ We call a category $\mathbf A $ \textit{trivially discrete}, if objects form a set,  and only morphisms in $\Mor(\mathbf A)$ are identity morphisms:
$$\Mor(\mathbf A)=\{1_{a}| a\in \Obj(\mathbf A)\}\simeq \Obj(\mathbf A).$$
Let $B$ be a topological space with an open cover $\{U_i\}_{i\in I},$ and ${B}^{\rm dis}$ be the corresponding trivially discrete category. Then, it is obvious that the set of trivially discrete categories, corresponding to $\{U_i\}_{i\in I},$ defines a categorical cover of  ${B}^{\rm dis}.$ 
}
\end{exmp}
\begin{exmp}\label{Ex:topcatcov2}
\rm{Suppose  $B:=M$ is a smooth manifold with an open cover $\{U_i\}_{i\in I}.$ Let $U_{ik}$ denotes the  intersection of $U_i$ and $U_k,$ 
$$U_{ik}:=U_i\cap U_k.$$
Let $\mathbf B:={\mathbb P} M$ be the path space category over $M;$ that is, (skipping the technicalities) 
\begin{equation}\nonumber
\begin{split}
&{\Obj}({\mathbb P}M)=M\\
&{\Mor}({\mathbb P}M)=\{C^1\hbox{-paths on $M$}\}/{\hbox{``some equivalence relations"}}.
\end{split}  
\end{equation}
Let, for each $i\in I,$ 
$${\mathbf U}_i:= {\mathbb P}U_i$$
be the path space category on the open set $U_i:$
\begin{equation}\nonumber
\begin{split}
&\Obj({{\mathbf U}_i})= U_i,\\
&\Mor({\mathbf U}_i)= \{C^1\hbox{-paths on $U_i$}\}/{\hbox{``some equivalence relations"}}.
\end{split}
\end{equation}
Clearly each ${\mathbf U}_i$ is a subcategory of ${\mathbb P} M.$ Observe that every path on $M$ has to locally lie on some $U_i.$ That means, for each path $\gamma$ on $M,$ we have $J\subset I,$  such that
$\gamma$ entirely lies on $\bigcup_{j\in J}U_j.$ Suppose $J=\{j_1, \cdots, j_k\},$ and $${\gamma}_{j_r}:=\gamma|_{U_{j_r}}=\hbox{restriction of} \hskip 0.15 cm \gamma \hskip 0.15 cm \hbox{to} \hskip 0.2 cm U_{j_r}.$$

So, $\gamma$ can be written as composition of a sequence of paths $(\gamma_{j_k}, \cdots, \gamma_{j_1})$:
$$\gamma=\gamma_{j_k}\circ\cdots \circ\gamma_{j_1} .$$
Note that $\gamma_{j_r}\in \Mor({\mathbf U}_{j_r}).$
In other words, according to the definition in \eqref{E:objmoruni}, $$\gamma\in \Mor\biggl(\bigcup_{j\in J} {\mathbf U}_j\biggr).$$
Also as  explained before each ${\mathbf U}_i$ is a groupoid. 
Hence, as a  consequence,  $\{{\mathbf U}_i\}_{i\in I}$ define a categorical cover of ${\mathbb P}M$,
$$\bigcup_{i\in I} {\mathbf U}_i={\mathbb P}M.\qquad \fbox{}$$

}
\end{exmp}
\begin{exmp}\label{Ex:topcatcov}
This is rather an example of what is \textit{not} an open categorical cover.
\rm{Let $B$ be a smooth topological space with an open cover $\{U_i\}_{i\in I}.$ Let $U_{ik}$ denotes the  intersection of $U_i$ and $U_k,$ 
$$U_{ik}:=U_i\cap U_k.$$
Let $\mathbf B$ be the path space grouopid  ${\mathbb P} B.$  Let $U_i^j,\hskip 0.2 cm i,j\in I,$ be the set of all morphisms whose sources are in $U_i$ and targets are in $U_j:$
\begin{equation}\label{E:uij}
U_i^j:=\{f\in \Mor({\mathbf B})| {\rm source}(f)\in U_i, {\rm target}(f)\in U_j\}\subset \Mor({\mathbf B}).
\end{equation}
Note that also,
$$U_i^j\cap U_k^l=U_{ik}^{jl},$$
where $U_{ik}^{jl}$ is the set of morphisms starting in $U_{ik},$ and terminating $U_{jl}.$  Since $\{U_i\}$ is an open cover of $B,$ every morphism in ${\mathbf B}$ has to start in some $U_i$ 
and terminate in some $U_j.$ In other words,
$$\bigcup_{i, j\in I}U_i^j=\Mor({\mathbf B}).$$
 Now let us define a category ${\mathbf U}_{i}^{k},$ for each $i, k\in I,$ as follows:
\begin{equation}\label{E:objmoruij}
\begin{split}
&\Obj({\mathbf U}_{i}^{j})={U}_i\cup {U}_j,\\
&\Mor({\mathbf U}_{i}^{j})={U}_i^j \cup \{1_p| p\in U_i\cup U_j\}.
\end{split}
\end{equation}
Thus each ${\mathbf U}_i^j$ is a subcategory of ${\mathbf U},$ and 
\begin{equation}\label{E:objmoruijcover}
\bigcup_{i, j\in I}{\mathbf U}_{i}^{j}={\mathbf B},
\end{equation}
But ${\mathbf U}_i^j$ is  not a groupoid. So according to our definition ${\mathbf U}_i^j$s are not open subcategories.
Hence, $\{{\mathbf U}_i^j\}_{i,j\in I}$ is not  an (open)categorical cover of ${\mathbf B}.$ \fbox{}
}
\end{exmp}


\section{$\mathbf {Cat}$-valued sheaves and functorial sections}\label{S:catsheaves}

In this section first we will introduce the notion of $\mathbf {Cat}$\textit{-valued sheaves}, and finally, our intention is to construct an example of $\mathbf {Cat}$ valued sheaf of functorial sections. 

\subsection{$\mathbf {Cat}$-valued sheaves}
Suppose   $\mathbf B$   is a topological groupoid. Let us consider the category ${\widetilde {\mathcal O}}(\mathbf B)$ of open subcategories:  

\begin{equation}\label{E:genobopen}
	\begin{split}
	  &{\rm Obj}\biggl({\widetilde {\mathcal O}}({\mathbf B})\biggr):=\{{\mathbf U}| {\mathbf U}\subset {\mathbf B}\}\\
	  &\hskip 2.5 cm =\hbox{set of all open subcategories of}\hskip 0.2 cm {\mathbf B},\\
	  &{\rm Hom}({\mathbf U}, {\mathbf V})=\{{\mbl}:{\mathbf U}\longrightarrow {\mathbf V}|{\mathbf U}, {\mathbf V}\subset {\mathbf B}\}.
\end{split}
\end{equation}
 Henceforth   ${\widetilde {\mathcal O}}(\mathbf B)$  will always denote the category of open subcategories of a topological groupoid $\mathbf B,$ as described in \eqref{E:genobopen}.
  
Recall the definition of a $\mathbf {Cat}$-valued presheaf, defined in \eqref{E:prshfcat}.  Let $\mathcal R$ be a $\mathbf {Cat}$-valued presheaf over ${\widetilde {\mathcal O}}(\mathbf B),$ 
\begin{equation}\label{E:preshtiob}
	{\mathcal R}:{\widetilde {\mathcal O}}{(\mathbf B) }^{\rm op}\to \rmc.
\end{equation}
Let ${\mathbf U}$ be an open subcategory of $\mathbf B,$ and $\{{\mathbf U}_{\alpha}\}_{\alpha\in I}$ be an open categorical cover of $\mathbf U$ (see \eqref{E:covobjmor}). Let
$${\mathbf i}_{\alpha}:{\mathbf U}_{\alpha}\hookrightarrow {\mathbf U}, \qquad \alpha\in I$$
be the inclusion functor.  Thus, ${\mathcal R}({\mathbf i}_{\alpha})$ defines a functor from the category ${\mathcal R}(\mathbf U)\in \Obj(\mathbf {Cat})$ to the category ${\mathcal R}({\mathbf U}_{\alpha})\in \Obj(\mathbf {Cat}):$
\begin{equation}\label{E:inclusheaf}
{\mathcal R}({\mathbf i}_{\alpha}): {\mathcal R}(\mathbf U)\longrightarrow {\mathcal R}({\mathbf U}_{\alpha}),
\end{equation}
for each $\alpha\in I.$ Similarly, if ${\mathbf U}_{\alpha}\cap {\mathbf U}_{\beta}:={\mathbf U}_{\alpha \beta}$ is non empty, then we have a pair of inclusion functors (we use same notation for both of them)
\begin{equation}\nonumber
\begin{split}
&{\mathbf i}_{\alpha \beta}: {\mathbf U}_{\alpha \beta}\hookrightarrow {\mathbf U}_{\alpha}, \\
&{\mathbf i}_{\alpha \beta}: {\mathbf U}_{\alpha \beta}\hookrightarrow {\mathbf U}_{\beta},
\end{split}
\end{equation}
and corresponding to them respectively we have  functors
\begin{equation}\label{E:ovlapinclu}
\begin{split}
{\mathcal R}({\mathbf i}_{\alpha \beta}): {\mathcal R}({\mathbf U}_{\alpha})\longrightarrow {\mathcal R}({\mathbf U}_{\alpha \beta}),\\
{\mathcal R}({\mathbf i}_{\alpha \beta}): {\mathcal R}({\mathbf U}_{\beta})\longrightarrow {\mathcal R}({\mathbf U}_{\alpha \beta}).
\end{split}
\end{equation}

We call $\mathcal R$ to be a $\mathbf {Cat}$-\textit{valued sheaf} over  ${\widetilde {\mathcal O}}(\mathbf B)$ provided following conditions are satisfied:
\begin{enumerate}
\item{{\bf {[Locality]}}
\begin{itemize}
\item[]{
(i)  If ${\mathbf \Psi}_1, {\mathbf \Psi}_2\in \Obj\biggl({\mathcal R}(\mathbf U)\biggr),$ such that,
\begin{equation}\label{E:objsheafcov}
{\mathcal R}({\mathbf i}_{\alpha})({\mathbf \Psi}_1)={\mathcal R}({\mathbf i}_{\alpha})({\mathbf \Psi}_2),
\end{equation}
for all $\alpha\in I,$ then,
\begin{equation}
{\mathbf \Psi}_1={\mathbf \Psi}_2.
\end{equation}
}
\item[]{
(ii) If ${\mathbf \Psi}_1, {\mathbf \Psi}_2, {\widetilde {\mathbf \Psi}}_1, {\widetilde {\mathbf \Psi}}_2 \in \Obj({\mathcal R}(\mathbf U))$ satisfy 
\begin{equation}\label{E:objsheafcovpair}
\begin{split}
&{\mathcal R}({\mathbf i}_{\alpha})({\mathbf \Psi}_1)={\mathcal R}({\mathbf i}_{\alpha})({\mathbf \Psi}_2),\\
&{\mathcal R}({\mathbf i}_{\alpha})({\widetilde {\mathbf \Psi}}_1)={\mathcal R}({\mathbf i}_{\alpha})({\widetilde {\mathbf \Psi}}_2),
\end{split}
\end{equation}
for all $\alpha\in I,$ and $\biggl({\mathbf \Psi}_1\xrightarrow {{\mathcal S}_1}{\widetilde {\mathbf \Psi}}_1\biggr), \biggl({\mathbf \Psi}_2\xrightarrow {{\mathcal S}_2} {\widetilde {\mathbf \Psi}}_2\biggr)\in \Mor\biggl({\mathcal R}(\mathbf U)\biggr)$ such that,
\begin{equation}\label{E:morsheafcov}
{\mathcal R}({\mathbf i}_{\alpha})({\mathcal S}_1)={\mathcal R}({\mathbf i}_{\alpha})({\mathcal S}_2),
\end{equation}
for all $\alpha\in I,$ then,
\begin{equation}
{\mathcal S}_1={\mathcal S}_2.
\end{equation}
}

\end{itemize}

}
\item{ {\bf {[Gluing]}}
\begin{itemize}
\item[]{
(i) If for each $\alpha\in I,$ a ${\mathbf \Psi}_\alpha\in \Obj\biggl({\mathcal R}({\mathbf U}_{\alpha})\biggr)$ is given, such that, for any non empty ${\mathbf U}_{\alpha \beta}$
\begin{equation}\label{E:objsheafover}
{\mathcal R}({\mathbf i}_{\alpha \beta})({\mathbf \Psi}_\alpha)={\mathcal R}({\mathbf i}_{\alpha \beta})({\mathbf \Psi}_\beta),
\end{equation}
for all $\alpha\in I,$ then there exists a ${\mathbf \Psi}\in \Obj({\mathbf U}),$ such that
\begin{equation}\nonumber
{\mathcal R}({\mathbf i}_{\alpha})({\mathbf \Psi})={\mathbf \Psi}_{\alpha}.
\end{equation}
}
\item[]{
(ii) If for each $\alpha\in I,$ a $\biggl({\mathbf \Psi}\xrightarrow{{\mathcal S}_\alpha} \widetilde {{\mathbf \Psi}}\biggr)\in \Mor\biggl({\mathcal R}({\mathbf U}_{\alpha})\biggr)$ is given, such that, for any non empty ${\mathbf U}_{\alpha \beta}$
\begin{equation}\label{E:objsheafoverpair}
\begin{split}
&{\mathcal R}({\mathbf i}_{\alpha \beta})({\mathbf \Psi}_\alpha)={\mathcal R}({\mathbf i}_{\alpha \beta})({\mathbf \Psi}_\beta),\\
&{\mathcal R}({\mathbf i}_{\alpha \beta})({\widetilde {\mathbf \Psi}}_\alpha)={\mathcal R}({\mathbf i}_{\alpha \beta})({\widetilde {\mathbf \Psi}}_\beta),
\end{split}
\end{equation}
and
\begin{equation}\label{E:ovmor}
{\mathcal R}({\mathbf i}_{\alpha \beta})({\mathcal S}_\alpha)={\mathcal R}({\mathbf i}_{\alpha \beta})({\mathcal S}_\beta).
\end{equation}
then there exists a $\biggl({\mathbf \Psi}\xrightarrow {\mathcal S}{\widetilde {\mathbf \Psi}}\biggr)\in \Mor({\mathbf U}),$ such that,
\begin{equation}\nonumber
{\mathcal R}({\mathbf i}_{\alpha})({\mathcal S})={\mathcal S}_{\alpha}.
\end{equation}
}

\end{itemize}

}

\end{enumerate}

\begin{exmp}\label{Ex:trivsheaf}
\rm{The above definition of a $\mathbf {Cat}$-valued sheaf is a  generalization of  the standard  definition of a sheaf.

 If we consider the trivially discrete (see Example~\ref{Ex:trivcov}   for the definition) topological category ${B}^{\rm dis}$ corresponding to a topological space $B,$ then a $\mathbf {Cat}$-valued sheaf over ${\mathcal O}({B}^{\rm dis})$ is simply a sheaf in the traditional sense.\fbox{}}
\end{exmp}

We end this section with a more interesting example of a $\mathbf {Cat}$-valued sheaf over ${\widetilde {\mathcal O}}({\mathbf B}).$

\subsection{${\mathbf {Cat}}$-valued sheaf of functorial sections}
Fix a category ${\mathbf C}.$ For any ${\mathbf U}\in \Obj\biggl({\widetilde {\mathcal O}}({\mathbf B})\biggr),$ let 
$${\mathcal F}(\mathbf U, \mathbf C):={\mathbf C}^{\mathbf U}$$
be the category of functors from $\mathbf U\longrightarrow \mathbf C;$ that is:
\begin{equation}\label{E:objmorsecsheaf}
\begin{split}
&\Obj({\mathbf C}^{\mathbf U})={\rm Fun}(\mathbf U, \mathbf C)\\
&\Mor({\mathbf C}^{\mathbf U})={\mathcal N}(\mathbf U, \mathbf C).
\end{split}
\end{equation}  
\begin{lemma}\label{L:escobjmor}
For any $\biggl({\mathbf U}\xrightarrow{\mathbf \Theta}{\mathbf V}\biggr)\in \Mor\biggl({\widetilde {\mathcal O}}(\mathbf B)\biggr),$ we have a functor
$${\mathcal R}(\mbl) : {\mathbf C}^{\mathbf V}\longrightarrow {\mathbf C}^{\mathbf U},$$ given as follows:
\begin{equation}\label{E:secobjmor}
\begin{split}
{\mathcal R}(\mbl) :  &\Obj({\mathbf C}^{\mathbf V})\longrightarrow \Obj({\mathbf C}^{\mathbf U}),\\
&{\mathbf \Psi}\mapsto {\mathbf \Psi}\mbl,\\
{\mathcal R}(\mbl) :  &\Mor({\mathbf C}^{\mathbf V})\longrightarrow \Mor({\mathbf C}^{\mathbf U}),\\
&{\mathcal S}\mapsto {\mathcal S}\mbl.
\end{split}
\end{equation}
\end{lemma}
\begin{proof}
Define a new category of categories ${\widehat {\mathcal C}}$ by attaching $\mathbf C$ with the original category ${\widetilde {\mathcal O}}(\mathbf B):$
\begin{equation}\nonumber
\begin{split}
&\Obj(\widehat{\mathcal C}):=\Obj\biggl({{\widetilde {\mathcal O}}(\mathbf B)}\biggr)\cup \{\mathbf C\},\\
&\Mor(\widehat{\mathcal C}):=\Mor\biggl({{\widetilde {\mathcal O}}(\mathbf B)}\biggr)\cup \biggl(\bigcup_{{\mathbf U}\subset \mathbf B}{\rm Fun}(\mathbf U, \mathbf C)\biggr)\cup \biggl(\bigcup_{{\mathbf U}\subset \mathbf B}{\rm Fun}(\mathbf C, \mathbf U)\biggr)\cup \biggl({\rm Fun}(\mathbf C, \mathbf C)\biggr).
\end{split}
\end{equation}
Then according to  Proposition~\ref{Pr:Hompresheaf}, for ${\mathbf C}\in \Obj({\widehat {\mathcal C}})$ we have the functor 
\begin{equation}
{\mathcal F}_{\mathbf C}:{\widehat{\mathcal C}}^{\rm op}\longrightarrow {\mathbf {Cat}}.
\end{equation}
In particular, if ${\mathbf U}, {\mathbf V}\in \Obj\bigl({\widetilde {\mathcal O}}(\mathbf B)\bigr)\subset \Obj(\widehat{\mathcal C})$ and $\biggl({\mathbf U}\xrightarrow{\mathbf \Theta}{\mathbf V}\biggr)$ then 
\begin{equation}\nonumber
\biggl({\mathbf U}\xrightarrow{\mathbf \Theta}{\mathbf V}\biggr)\mapsto  \biggl({\mathbf C}^{\mathbf V}\xrightarrow{{\mathcal F}_{\mathbf C}(\mbl)} {\mathbf C}^{\mathbf U}\biggr),
\end{equation}
since ${\mathcal F}_{\mathbf C}(\mathbf U)={\mathcal F}(\mathbf U, \mathbf C)={\mathbf C}^{\mathbf U}.$
It is immediate that ${\mathcal F}_{\mathbf C}$ restricted to ${\widetilde {\mathcal O}}(\mathbf B)$ is in fact the $\mathcal R$ in the statement of the lemma;
$${\mathcal F}_{\mathbf C}|_{{\widetilde {\mathcal O}}(\mathbf B)}=\mathcal R.$$

\end{proof}

Therefore we have a $\mathbf {Cat}$-valued presheaf $\mathcal R:$
\begin{equation}\label{E:presheafsec}
\begin{split}
{\mathcal R} : &{\widetilde {\mathcal O}}({\mathbf B})^{\rm op}\to {\mathbf {Cat}},\\
&\Obj\biggl({\widetilde {\mathcal O}}({\mathbf B})\biggr)\to \Obj({\mathbf {Cat}}),\\
&{\mathbf U}\mapsto {\mathbf C}^{\mathbf U},\\
&\Mor\biggl({\widetilde {\mathcal O}}({\mathbf B})\biggr)\to \Mor(\mathbf {Cat}),\\
&\biggl({\mathbf U}\xrightarrow{\mbl}{\mathbf V}\biggr)\mapsto \biggl({\mathbf C}^{\mathbf V}\xrightarrow{{{\mathcal R}(\mbl)}} {\mathbf C}^{\mathbf U}\biggr).
\end{split}
\end{equation}
Rest of this section will be devoted to prove  (Theorem~\ref{Th:funcsheaf}) that the $\mathbf {Cat}$-valued presheaf, defined above, is in fact a  $\mathbf {Cat}$-valued sheaf.

\begin{theorem}\label{Th:funcsheaf}
Let $\mathbf B$ be a topological groupoid. Let $\widetilde {\mathcal O}(\mathbf B)$ be the category of open subcategories.  Let $\mathcal R$ be as defined in Lemma~\ref{L:escobjmor}. Then $\mathcal R$ is a $\mathbf {Cat}$-valued sheaf over $\widetilde {\mathcal O}(\mathbf B).$
\end{theorem}
Note that if $\{{\mathbf U}_{\alpha}\}_{\alpha}$ is an open categorical cover of $\mathbf B,$ then  Lemma~\ref{L:uniintprop} ensures,
\begin{equation}\label{E:uinitprop}
\begin{split}
&{\mathbf U}_{\alpha \beta}:={\mathbf U}_{\alpha}\cap {\mathbf U}_{\beta}\in \Obj\bigl({\widetilde {\mathcal O}}(\mathbf B)\bigr),\\
&{\mathbf U}_{\alpha}\cup{\mathbf U}_{\beta}\in \Obj\bigl({\widetilde {\mathcal O}}(\mathbf B)\bigr).
\end{split}
\end{equation}

Before we carry out  the proof of this theorem,  we will work in a simpler situation to  get an insight into  the overall methodology of the proof.

We have to verify that $\mathcal R$ satisfies {\bf {[Locality]}} and  {\bf {[Gluing]}} conditions.  The crucial issue here is that, if $\{{\mathbf U}_{\alpha}\}$ is a categorical cover of $\mathbf U$, then in general $\bigcup_{\alpha}\Mor({\mathbf U}_{\alpha})$ is a subset of $\Mor\biggl(\bigcup_{\alpha}{\mathbf U}_{\alpha}\biggr)=\Mor(\mathbf U).$ In other words, there may exist a morphism in $\Mor({\mathbf U}),$ which does not belong to any $\Mor({\mathbf U}_{\alpha}).$
Now, if ${\mathbf \Psi}_1, {\mathbf \Psi}_2 \in \Obj({\mathcal R}(\mathbf U))=\Obj({\mathbf C}^{\mathbf U})={\rm Fun}(\mathbf U, \mathbf C),$ such that ${\mathbf \Psi}_1, {\mathbf \Psi}_2$ coincide in each ${\mathbf U}_{\alpha},$ 
then it is not obvious that  ${\mathbf \Psi}_1={\mathbf \Psi}_2.$ Because, there is a possibility that ${\mathbf \Psi}_1, {\mathbf \Psi}_2$ are different on some 
\begin{equation}\label{E:notin}
\begin{split}
&f\in \Mor\biggl(\bigcup_{\alpha}{\mathbf U}_{\alpha}\biggr)=\Mor(\mathbf U),\\ 
&f\notin \bigcup_{\alpha}\Mor({\mathbf U}_{\alpha}).
\end{split}
\end{equation}
 But we will see, remarkably, that does not actually happen! And, even for  an $f$ as in \eqref{E:notin},  ${\mathbf \Psi}_1(f), {\mathbf \Psi}_2(f)$ are completely determined by the local datum; i.e. if \eqref{E:objsheafcov}
holds for ${\mathbf \Psi}_1, {\mathbf \Psi}_2,$ then ${\mathbf \Psi}_1(f)={\mathbf \Psi}_2(f),$ for any $f\in \Mor(\mathbf U)$ including  one like in \eqref{E:notin}. Similarly,  other conditions of {\bf {[Locality]}} and  {\bf {[Gluing]}} also hold. To see this, at first  we deal with a test case, where an open subcategory $\mathbf U$ of ${\mathbf B}$ is covered by only two (open) subcategories $\{{\mathbf U}_1, {\mathbf U}_2\}:$
\begin{equation}\label{E:cover2}
{\mathbf U}_1\cup {\mathbf U}_2=\mathbf U.
\end{equation}
Before we proceed further, let us simplify our notations. If $\mathbf U$ is a subcategory of $\mathbf V,$ and ${\mathbf i}:{\mathbf U}\hookrightarrow \mathbf V$ is the inclusion functor, then for  ${\mathbf \Psi}\in \Obj({\mathbf C}^{\mathbf V}), {\mathcal S}\in \Mor({\mathbf C}^{\mathbf V}),$ we will respectively denote,
\begin{equation}\label{E:rest}
\begin{split}
\biggl({\mathcal R}(\mathbf i)\biggr)(\mathbf \Psi):={\mathbf \Psi}|_{\mathbf U}\in \Obj({\mathbf C
}^{\mathbf U}),\\
\biggl({\mathcal R}(\mathbf i)\biggr)(\mathcal S):={\mathcal S}|_{\mathbf U}\in \Mor({\mathbf C
}^{\mathbf U}).
\end{split}
\end{equation} 
\begin{prop}\label{pr:cover2sheaf}
Let $\widetilde {\mathcal O}(\mathbf B)$ be the category of open subcategories.  Let $\mathcal R$ be as defined in Lemma~\ref{L:escobjmor}. Let ${\mathbf U}$ be an open subcategory of $\widetilde {\mathcal O}(\mathbf B)$ with  categorical cover $\{{\mathbf U}_{\alpha}\}_{\alpha\in \{1, 2\}}.$ Then on $\mathbf U$ {\bf {[Locality]}} and  {\bf {[Gluing]}} conditions, listed between \eqref{E:objsheafcov}--\eqref{E:ovmor}, hold. 
\end{prop}
\begin{proof} 
First we verify {\bf {[Locality]}} conditions.

Let ${\mathbf \Psi}_1, {\mathbf \Psi}_2\in \Obj({\mathbf C}^{\mathbf U})={\rm Fun}(\mathbf U, \mathbf C),$ such that,
$${\mathbf \Psi}_1|_{{\mathbf U}_{\alpha}}= {\mathbf \Psi}_2|_{{\mathbf U}_{\alpha}}, \qquad\alpha\in \{1, 2\}.$$
Recalling the definition of categorical union given in \eqref{E:objmoruni} we see:
$\mbeta_1=\mbeta_2$ on objects. Also, $\mbeta_1={\mathbf \Psi}_2$ on $\Mor({\mathbf U}_1)\cup \Mor({\mathbf U}_2)\subset \Mor({\mathbf U}).$ Now, consider an arbitrary $g\in \Mor(\mathbf U),$ that means there exist $g_2\in \Mor({\mathbf U}_2)$ and $g_1\in \Mor({\mathbf U}_1),$ such that
$$g=g_2\circ g_1.$$
[ To be precise, ${\tilde g}_1\in \Mor({\mathbf U}_1), {\tilde g}_2\in \Mor({\mathbf U}_2),$ and ${\tilde g}_1\circ {\tilde g}_2=g$ is also an alternate possibility. But, this  case can be dealt in exactly same fashion as the other. And, we will not explicitly consider it.]
Thus, 
\begin{equation}\label{E:equobj}
\begin{split}
\mbeta_1(g)&=\mbeta_1(g_2)\circ \mbeta_1(g_1) \hskip 0.4 cm [\hskip 0.05 cm \hbox{since}, \mbeta_1 \hbox{is a functor}]\\
&= \mbeta_2(g_2)\circ \mbeta_2(g_1) \hskip 0.2 cm [\hskip 0.05 cm \hbox{since},  g_2\in \Mor({\mathbf U}_2), g_1\in \Mor({\mathbf U}_1),\\ 
&\hskip 4 cm \hbox{and}  \hskip 0.15 cm {\mathbf \Psi}_1|_{{\mathbf U}_{\alpha}}= {\mathbf \Psi}_2|_{{\mathbf U}_{\alpha}}, \alpha\in \{1, 2\}]\\
&= \mbeta_2(g_2\circ g_1) \hskip 0.4 cm [\hskip 0.05 cm \hbox{since}, \mbeta_2\hskip 0.15 cm \hbox{is a functor}]\\
&=\mbeta_2(g).
\end{split}
\end{equation}
Thus $\mbeta_1=\mbeta_2$ on $\mathbf U.$

Similarly, if $\bigl({\mathcal S}_1:\mbeta_1\Longrightarrow \mbeta_2\bigr), \bigl({\mathcal S}_2:\mbeta_1\Longrightarrow \mbeta_2\bigr)\in \Mor({\mathbf C}^{\mathbf U})={\mathcal N}(\mathbf U, \mathbf C),$ such that,
$${\mathcal S}_1|_{{\mathbf U}_{\alpha}}= {\mathcal S}_2|_{{\mathbf U}_{\alpha}}, \qquad\alpha\in \{1, 2\},$$
then for any $(a\xrightarrow{f}b)\in \Mor({\mathbf U}_1)\cup \Mor({\mathbf U}_2),$ we have the commuting diagram,
\begin{equation} \label{D:morsheaf}
\xymatrix{
         \ar[d]^{{\mathcal S}(a)} \mbeta_1(a)     \ar[rr]^-{\mbeta_1(f)} & &\mbeta_1(b) \ar[d]_{{\mathcal S}(b)} \\
\mbeta_2(a) \ar[rr]_-{\mbeta_2(f)}& & \mbeta_2(b) 
},
\end{equation}
where we denote ${\mathcal S}:={\mathcal S}_1|_{{\mathbf U}_{\alpha}}= {\mathcal S}_2|_{{\mathbf U}_{\alpha}}, \alpha\in \{1, 2\}.$ Now suppose $(a\xrightarrow{g}c)\in \Mor(\mathbf U)$ is an arbitrary morphism in $\Mor(\mathbf U).$ That means we have $(a\xrightarrow{g_1}b)\in \Mor({\mathbf U}_1), (b\xrightarrow{g_2}c)\in \Mor({\mathbf U}_2),$ such that
$$g=g_2\circ g_1.$$
$\mathcal S$ extends to $(a\xrightarrow{g}c)\in \Mor(\mathbf U)$ as follows. We have commuting diagrams,
\begin{equation} \label{D:morsheafside1}
\xymatrix{
         \ar[d]^{{\mathcal S}(a)} \mbeta_1(a)     \ar[rr]^-{\mbeta_1(g_1)} & &\mbeta_1(b) \ar[d]_{{\mathcal S}(b)} \\
\mbeta_2(a) \ar[rr]_-{\mbeta_2(g_1)}& & \mbeta_2(b) 
},
\end{equation}
and,

\begin{equation} \label{D:morsheafside2}
\xymatrix{
         \ar[d]^{{\mathcal S}(b)} \mbeta_1(b)     \ar[rr]^-{\mbeta_1(g_2)} & &\mbeta_1(c) \ar[d]_{{\mathcal S}(c)} \\
\mbeta_2(b) \ar[rr]_-{\mbeta_2(g_2)}& & \mbeta_2(c) 
}.
\end{equation}
Since $\mbeta_2, \mbeta_1$ are functors, composing above commuting diagrams we get the commuting diagram,
\begin{equation} \label{D:morsheafside3}
\xymatrix{
         \ar[d]^{{\mathcal S}(a)} \mbeta_1(a)     \ar[rr]^-{\mbeta_1(g)} & &\mbeta_1(c) \ar[d]_{{\mathcal S}(c)} \\
\mbeta_2(a) \ar[rr]_-{\mbeta_2(g)}& & \mbeta_2(c) 
} \qquad[\hskip 0.06 cm \hbox{putting}\hskip 0.15 cm g_2\circ g_1=g ].
\end{equation}
Thus ${\mathcal S}={\mathcal S}_1={\mathcal S}_2$ for the category $\mathbf U.$

Our next task is to verify {\bf {[Gluing]}} conditions.

For that, suppose $\mbeta_{\alpha}\in \Obj({\mathbf C}^{{\mathbf U}_{\alpha}})={\rm Fun}({\mathbf U}_{\alpha}, {\mathbf C})$ given for $\alpha\in \{1, 2\},$ such that:
$\mbeta_{1}|_{{\mathbf U}_{1 2}}=\mbeta_{2}|_{{\mathbf U}_{1 2}}.$
Then we define $(\mbeta:{\mathbf U}\longrightarrow {\mathbf C}) \in {\rm Fun}({\mathbf U}, {\mathbf C})=\Obj({\mathbf C}^{\mathbf U})$ as follows:
\begin{eqnarray}
&&\Obj({\mathbf U})\to \Obj({\mathbf C})\nonumber\\
&&\Obj({\mathbf U}_1)\cup \Obj({\mathbf U}_2)\to \Obj({\mathbf C})\nonumber\\
&&\hskip 3 cm a\mapsto \mbeta_{\alpha}(a), \qquad \hbox{if}\hskip 0.15 cm  a\in {\mathbf U}_{\alpha}, \alpha\in\{1, 2\},\label{E:gluobj}
\end{eqnarray}
and,
\begin{eqnarray}
&&\Mor({\mathbf U})\to \Mor({\mathbf C})\nonumber\\
&&{\mathbf {Gen}}\biggl(\Mor({\mathbf U}_1)\cup \Mor({\mathbf U}_2)\biggr)\to \Mor({\mathbf C}),\nonumber\\
&&\hskip 3.5 cm f_2\circ f_1\mapsto \mbeta_2(f_2)\circ \mbeta_1(f_1), \nonumber\\
&&\hskip 3.8 cm f_2\in \Mor({\mathbf U}_2), f_1\in \Mor({\mathbf U}_1).\label{E:glumor}
\end{eqnarray}
Since $\mbeta_{1}|_{{\mathbf U}_{1 2}}=\mbeta_{2}|_{{\mathbf U}_{1 2}},$ \eqref{E:gluobj} is well defined. For the same reason, right hand side of \eqref{E:glumor} makes sense [$t(f_1)=s(f_2)\in \Obj({\mathbf U}_{1 2}).$ So, $s(\mbeta_2(f_2))= t(\mbeta_1(f_1)).$]. But we have to check if \eqref{E:glumor} is well defined, that is, if $f'_2\circ f'_1=f_2\circ f_1,$ for some other $f'_2\in \Mor({\mathbf U}_2), f'_1\in \Mor({\mathbf U}_1),$ then we should have
$$\mbeta_2(f_2)\circ \mbeta_1(f_1)= \mbeta_2(f'_2)\circ \mbeta_1(f'_1).$$
Let $$s(f_2)=b=t(f_1), s(f_1)=a=s(f'_1), t(f'_2)=c=t(f_2), s(f'_2)=b'=t(f'_1),$$
so $b, b'\in \Obj({\mathbf U}_{1 2}).$ 

By Lemma~\ref{L:property}  we have an isomorphism $b\xrightarrow{g_0}b'\in \Mor(\mathbf U)$ such that:
\begin{equation}\nonumber
\begin{split}
&g_0\circ f_1=f'_1,\\
&f_2\circ g_0^{-1}=f'_2.
\end{split}
\end{equation}
In fact,  $$b\xrightarrow{g_0}b'\in \Mor({\mathbf U}_{1 2}).$$
Then 
\begin{equation}
\begin{split}
 &\mbeta_2(f'_2)\circ \mbeta_1(f'_1)\\
 &= \mbeta_2(f_2 \circ g_0^{-1})\circ \mbeta_1(g_0\circ f_1)\\
 &= \mbeta_2(f_2) \circ \mbeta_2(g_0^{-1})\circ \mbeta_1(g_0)\circ \mbeta_1(f_1)\\
 &=\mbeta_2(f_2)\circ \mbeta_1(f_1) \hskip 0.4 cm [\hskip 0.06 cm \hbox{since}\hskip 0.15 cm g_0\in \Mor({\mathbf U}_{1 2})\hskip 0.15 cm \hbox{and}\hskip 0.15 cm \mbeta_{1}|_{{\mathbf U}_{1 2}}=\mbeta_{2}|_{{\mathbf U}_{1 2}}] .
 \end{split}
 \end{equation}
 
Hence, \eqref{E:glumor} is well defined, and we have $(\mbeta:{\mathbf U}\longrightarrow {\mathbf C}) \in {\rm Fun}({\mathbf U}, {\mathbf C})=\Obj({\mathbf C}^{\mathbf U})$ satisfying
\begin{equation}\nonumber
\mbeta |_{{\mathbf U}_{\alpha}}=\mbeta_{\alpha}, \qquad \alpha=\{1, 2\}
\end{equation}
Next, suppose $\bigl({\mathcal S}_{\alpha}:{\mbeta}_{\alpha}\Longrightarrow \mbeta'_{\alpha}\bigr)\in \Mor({\mathbf C}^{{\mathbf U}_{\alpha}})={\mathcal N}({\mathbf U}_{\alpha}, {\mathbf C})$ given for $\alpha\in \{1, 2\},$ such that,
$\mbeta_{1}, \mbeta_{2}$ glue to form a $\mbeta\in {\rm Fun}(\mathbf U, \mathbf C),$ similarly, $\mbeta'_{1}, \mbeta'_{2}$ glue to form a $\mbeta'\in {\rm Fun}(\mathbf U, \mathbf C),$ [as described in the previous part of the proof] and ${\mathcal S}_{\alpha}$ satisfy 
\begin{equation}\label{E:natov}
{\mathcal S}_{1}|_{{\mathbf U}_{1 2}}={\mathcal S}_{2}|_{{\mathbf U}_{1 2}}.
\end{equation}
We define ${\mathcal S}\in {\mathcal N}({\mathbf U}, {\mathbf C})$ by:
\begin{equation}\label{E:glonat}
{\mathcal S}(a)={\mathcal S}_{\alpha}(a), \hskip 0.4 cm \hbox{for}\hskip 0.15 cm a\in \Obj({\mathbf U}_{\alpha}), \alpha\in \{1, 2\}.
\end{equation}
Above equation makes sense because of \eqref{E:natov}. We have to ensure that \eqref{E:glonat} is a well defined natural transformation between $\mbeta$ and $\mbeta'.$
It is obvious that for any $a\xrightarrow{f}b\in \Mor({\mathbf U}_1)\cup \Mor({\mathbf U}_2)\subset \Mor(\mathbf U)=\Mor({\mathbf U}_1\cup {\mathbf U}_2),$ we have the commuting diagram:
\begin{equation} \label{D:morsheafglu}
\xymatrix{
         \ar[d]^{{\mathcal S}_{\alpha}(a)} \mbeta(a)     \ar[rr]^-{\mbeta(f)} & &\mbeta(b) \ar[d]_{{\mathcal S}_{\alpha}(b)} \\
\mbeta'(a) \ar[rr]_-{\mbeta'(f)}& & \mbeta'(b) 
},
\end{equation}
where $a\xrightarrow{f}b\in \Mor({\mathbf U}_\alpha), \alpha=\{1, 2\}.$  We have to  verify that ${\mathcal S}$ is well defined (as a natural transformation) for all $(a\xrightarrow{g}c)\in \Mor(\mathbf U).$ Again, by the previous argument, we have $g=g_2\circ g_1,$ for some $(a\xrightarrow{g_1}b)\in \Mor({\mathbf U}_1), (b\xrightarrow{g_2}c)\in \Mor({\mathbf U}_2).$ We have a pair of commuting diagrams respectively in ${\mathbf U}_1$ and ${\mathbf U}_2,$
\begin{equation} \label{D:morsheafglucom1}
\xymatrix{
         \ar[d]^{{\mathcal S}_{1}(a)} \mbeta(a)     \ar[rr]^-{\mbeta(g_1)} & &\mbeta(b) \ar[d]_{{\mathcal S}_{1}(b)} \\
\mbeta'(a) \ar[rr]_-{\mbeta'(g_1)}& & \mbeta'(b) 
},
\end{equation}
and,

\begin{equation} \label{D:morsheafglucom2}
\xymatrix{
         \ar[d]^{{\mathcal S}_{2}(b)} \mbeta(b)     \ar[rr]^-{\mbeta(g_2)} & &\mbeta(c) \ar[d]_{{\mathcal S}_{2}(c)} \\
\mbeta'(b) \ar[rr]_-{\mbeta'(g_2)}& & \mbeta'(c) 
}.
\end{equation}
Since $b\in \Obj({\mathbf U}_{1 2}),$ and ${\mathcal S}_{1}|_{{\mathbf U}_{1 2}}={\mathcal S}_{2}|_{{\mathbf U}_{1 2}},$ we have a  commuting diagram
\begin{equation} \label{D:morsheafglucom}
\xymatrix{
         \ar[d]^{{\mathcal S}_{1}(a)} \mbeta(b)     \ar[rr]^-{\mbeta(g)} & &\mbeta(c) \ar[d]_{{\mathcal S}_{2}(c)} \\
\mbeta'(b) \ar[rr]_-{\mbeta'(g)}& & \mbeta'(c) 
} \qquad [\hskip 0.06 cm \hbox{putting} \hskip 0.15 cm g_2\circ g_1=g].
\end{equation}
So we have a well defined natural transformation $(\mathcal S:\mbeta\Longrightarrow \mbeta')\in {\mathcal N}(\mathbf U, \mathbf C)$ satisfying:
$${\mathcal S}|_{{\mathbf U}_{\alpha}}={\mathcal S}_{\alpha}, \qquad \alpha\in\{1, 2\}.$$

In conclusion, we have proven that,  on ${\mathbf U}={\mathbf U}_1\cup {\mathbf U}_2, \mathcal R$ satisfies {\bf {[Locality]}} and  {\bf {[Gluing]}} conditions.

\end{proof}

Now we turn to the proof of Theorem~\ref{Th:funcsheaf}. The proof would be rather long, and we proceed by one step at a time. The verification of {\bf {[Locality]}} is not very difficult. The main hardship is to ensure  that the {\bf {[Gluing]}} conditions are satisfied. Let us first explain what we are trying to achieve here, and give a brief description of the strategy of our proof.

Suppose $\{{\mathbf U}_{\alpha}\}_{\alpha\in I}$ is a categorical cover of an open subcategory $\mathbf U\subset {\mathbf B}.$ Let $f\in \Mor({\mathbf U})=\Mor(\cup_{\alpha\in I}{\mathbf U}_{\alpha})$ be an arbitrary morphism in $\Mor(\mathbf U).$
That means, there exists $J:=\{j_1, \cdots, j_m\}\subset I$ such that
\begin{equation}\label{E:comcoverj}
f=f_{j_m}\circ\cdots \circ f_{j_1},
\end{equation}
where $f_{j_r}\in \Mor({\mathbf U}_{j_r}), r\in\{1, \cdots, m\}.$ Note that ${\mathbf U}_{j_r}\cap {\mathbf U}_{j_{r-1}}, r\in \{2,\cdots, m\},$ is always non empty, because 
$$s(f_{j_r})=t(f_{j_{r-1}}).$$ 
If we are given a $\mbeta_{\alpha}\in \Obj({\mathbf C}^{{\mathbf U}_{\alpha}})={\rm Fun}({\mathbf U}_{\alpha}, \mathbf C)$ for each $\alpha\in I,$  satisfying 
$$\mbeta_{\alpha}|_{{\mathbf U}_{\alpha \beta}}=\mbeta_{\beta}|_{{\mathbf U}_{\alpha \beta}}\qquad \hbox{for any non-empty}\hskip 0.15 cm {\mathbf U}_{\alpha \beta},$$
we intend   to find a $\mbeta\in \Obj({\mathbf C}^{\mathbf U})={\rm Fun}(\mathbf U, \mathbf C),$ such that
$$\mbeta|_{{\mathbf U}_{\alpha}}=\mbeta_{\alpha}.$$
We define $\mbeta$ as follows
\begin{equation}\label{E:globalpsij}
\mbeta(f):=\mbeta_{j_m}(f_{j_m})\circ\cdots\circ \mbeta_{j_1}(f_{j_1}).
\end{equation}
[ Here and afterwards we will focus  on the gluing (respectively locality) condition only for morphisms.  For the objects it is obvious due to the first part of \eqref{E:objmoruni}.]

Now, suppose for some other $J':=\{j'_1, \cdots, j'_n\}\subset I,$ we have
\begin{equation}\label{E:comcoverj'}
f=f'_{j'_n}\circ\cdots \circ f'_{j'_1},
\end{equation}
where $f'_{j'_k}\in \Mor({\mathbf U}_{j'_k}), k\in\{1, \cdots, n\}.$ Then $\mbeta$ to be well defined, we should have
\begin{equation}\label{E:globalpsijij'}
\mbeta_{j_m}(f_{j_m})\circ\cdots\circ \mbeta_{j_1}(f_{j_1})=\mbeta_{j'_n}(f'_{j'_n})\circ\cdots\circ \mbeta_{j'_1}(f'_{j'_1}).
\end{equation}
Our primary concern is to prove \eqref{E:globalpsijij'}. Rest of the proof is straightforward. 
\subsubsection*{{\bf {Strategy of the proof}}}
The proof of \eqref{E:globalpsijij'} will be carried out by the method of  induction. First we will prove that, if $J=\{j_1\}$ and $J'=\{j'_2, j'_1\},$ that is,  $f\in \Mor({\mathbf U}_{j_1})$ and $f=f'_{j'_2}\circ f'_{j'_1},$ then
$$\mbeta_{j_1}(f)=\mbeta_{j'_2}(f'_{j'_2})\circ\mbeta_{j'_1}(f'_{j'_1}).$$
Next, we make the induction assumption:

\subsubsection*{\bf {[Ind assum 1]:}} ``If $J=\{j_1\}$ and $J'=\{j'_1, \cdots, j'_k\},$ that is,  $f\in \Mor({\mathbf U}_{j_1})$ and $f=f'_{j'_k}\circ\cdots\circ f'_{j'_1},$ then
$$\mbeta_{j_1}(f)=\mbeta_{j'_k}(f'_{j'_k})\circ\cdots\circ\mbeta_{j'_1}(f'_{j'_1}).\hbox{''}$$

Then we show that  it is true for $J=\{j_1\}$ and $J'=\{j'_{1}, \cdots, j'_{k+1}\}.$ Thus, it holds for $J=\{j_1\}$ and $J'=\{j'_1, \cdots, j'_n\}.$ We proceed with a second   induction assumption.

\subsubsection*{\bf {[Ind assum 2]:}} ``If $J=\{j_1,\cdots, j_{k}\}$ and $J'=\{j'_1, \cdots, j'_n\},$ that is,  $f_{j_l}\in \Mor({\mathbf U}_{j_l})$ and $f_{j_k}\circ\cdots\circ f_{j_1}=f'_{j'_n}\circ\cdots\circ f'_{j'_1},$ then
$$\mbeta_{j_k}(f_{j_k})\circ \cdots \circ \mbeta_{j_1}(f_{j_1})=\mbeta_{j'_n}(f'_{j'_n})\circ\cdots\circ\mbeta_{j'_1}(f'_{j'_1}).\hbox{''}$$
Then we show that it holds for $J=\{j_1,\cdots, j_{k+1}\}$ and $J'=\{j'_1, \cdots, j'_n\}.$ Thus, we can conclude that it holds for $J=\{j_1,\cdots, j_{m}\}$ and $J'=\{j'_1, \cdots, j'_n\}.$

And, that would complete the proof of \eqref{E:globalpsijij'}.
\subsubsection*{{\bf {Proof of \eqref{E:globalpsijij'}}}} We implement our strategy described above.
\begin{prop}\label{pr:glusheaf2and1}
With notations  and conventions as above, let ${\mathbf U}$ be a subcategory of $\mathbf B,$ with an open categorical cover $\{{\mathbf U}_{\alpha}\}_{\alpha\in I}.$ Let $j_1, j'_1, j'_2 \in I,$ and 
$$f=f'_{j'_2}\circ f'_{j'_1}\in \Mor({{\mathbf U}_{j_1}})\cap \Mor({\mathbf U}_{j'_1}\cup {\mathbf U}_{j'_2}),$$
where $f\in \Mor({\mathbf U}_{j_1}), f'_{j'_2}\in \Mor({\mathbf U}_{j'_2}), f'_{j'_1}\in \Mor({\mathbf U}_{j'_1}).$ Then,
\begin{equation}\label{E:glusheafkand1}
\mbeta_{j_1}(f)=\mbeta_{j'_2}(f'_{j'_2})\circ \mbeta_{j'_1}(f'_{j'_1}).
\end{equation}
\end{prop}
\begin{proof}
First we observe that, both ${\mathbf U}_{j_1}\cap {\mathbf U}_{j'_1}$ and ${\mathbf U}_{j_1}\cap {\mathbf U}_{j'_2}$ are non empty [since, $s(f)=s(f'_{j'_1}), t(f)=t(f'_{j'_2}),$ and $s(f), t(f)\in \Obj({\mathbf U}_{j_1})$].  Now, using (iii) of Proposition~\ref{pr:boolean}, we write
$$f \in   \Mor\biggl({\mathbf U}_{j_1}\cap \bigl({\mathbf U}_{j'_1} \cup {\mathbf U}_{j'_2} \bigr)\biggr)=\Mor\biggl(\bigl({\mathbf U}_{j_1}\cap {\mathbf U}_{j'_1}\bigr)\cup \bigl({\mathbf U}_{j_1}\cap {\mathbf U}_{j'_2}\bigr)\biggr).$$
That means, there exists $f_2\in \Mor\bigl({\mathbf U}_{j_1}\cap {\mathbf U}_{j'_2}\bigr), f_1\in \Mor\bigl({\mathbf U}_{j_1}\cap {\mathbf U}_{j'_1}\bigr),$ such that
$$f=f_2\circ f_1.$$
By \eqref{E:samecompo}--\eqref{E:compocon}, we have an isomorphism $g_0\in \Mor({\mathbf U}_{j'_1}\cap {\mathbf U}_{j'_2}),$ such that
\begin{equation}\nonumber
\begin{split}
f'_{j_2}\circ g_0=f_2,\\
{{g_0}^{-1}}\circ f'_{j_1}=f_1.
\end{split}
\end{equation}
So,
\begin{equation}\label{E:2t01fin}
\begin{split}
{\mbeta}_{j_1}(f)&={\mbeta}_{j_1}(f_2\circ f_1)\\
&= {\mbeta}_{j_1}(f_2)\circ {\mbeta}_{j_1}(f_1) [\hskip 0.05 cm \hbox{since}, {\mbeta}_{j_1}\hbox{ is a functor, and}\hskip 0.15 cm  f_2, f_1\in \Mor({\mathbf U}_{j_1})]\\
&= {\mbeta}_{j'_2}(f_2)\circ {\mbeta}_{j'_1}(f_1)\\
& [\hskip 0.05 cm \hbox{since},  {\mbeta}_{j_1}|_{{\mathbf U}_{j_1 j'_i}}={\mbeta}_{j'_i}|_{{\mathbf U}_{j_1 j'_i}}, i\in\{1, 2\}, \hskip 0.15 cm \hbox{and}\hskip 0.15 cm f_2 \in \Mor({\mathbf U}_{j'_2}), f_1\in \Mor({\mathbf U}_{j'_1})]\\
&={\mbeta}_{j'_2}(f'_{j_2}\circ g_0)\circ {\mbeta}_{j'_1}( {{g_0}^{-1}}\circ f'_{j_1})\\
&={\mbeta}_{j'_2}(f'_{j'_2})\circ {\mbeta}_{j'_1}(f'_{j'_1}) [\hskip 0.05 cm \hbox{since},  {\mbeta}_{j'_1}|_{{\mathbf U}_{j'_1 j'_2}}={\mbeta}_{j'_2}|_{{\mathbf U}_{j'_1 j'_2}}]
\end{split}
\end{equation}
This proves the proposition.
\end{proof}

We assume {\bf {[Ind assum 1]}}.

\begin{prop}\label{pr:glusheafk+1and1}
Suppose {\bf {[Ind assum 1]}} is true.  Let $j_1\in I, {\widetilde J}'=\{j'_1,\cdots, j'_{k+1}\} \subset I,$ and 
$$f=f'_{j'_{k+1}}\circ\cdots\circ f'_{j'_1}\in \Mor({{\mathbf U}_{j_1}})\bigcap \Mor\biggl(\bigcup_{j'_i\in {\widetilde J}'}{\mathbf U}_{j'_1}\biggr),$$
where $f\in \Mor({\mathbf U}_{j_1}), f'_{j'_i}\in \Mor({\mathbf U}_{j'_i}), i\in\{1,\cdots, k+1\}.$ Then,
\begin{equation}\label{E:glusheaf2and1}
\mbeta_{j_1}(f)=\mbeta_{j'_{k+1}}(f'_{j'_{k+1}})\circ \cdots\circ\mbeta_{j'_1}(f'_{j'_1}).
\end{equation}

\end{prop}
\begin{proof}
In spirit, proof is similar to that of Proposition~\ref{pr:glusheaf2and1}. 
Let us denote $${\mathbf U}_{j'_1}\cup\cdots\cup{\mathbf U}_{j'_k} \stackrel{\small{{\rm {notation}}}}{=}{\mathbf U}(j'_1,\cdots, j'_k).$$
Again  ${\mathbf U}_{j_1}\cap {\mathbf U}({j'_1\cdots j'_k})$ and ${\mathbf U}_{j_1}\cap {\mathbf U}_{j'_{k+1}}$ are non empty [since, $s(f)=s(f'_{j'_1}), t(f)=t(f'_{j'_{k+1}}),$ and $s(f), t(f)\in \Obj({\mathbf U}_{j_1})$].
Then we have,
$$f \in   {\mathbf U}_{j_1}\bigcap \biggl({\mathbf U}({j'_1\cdots j'_k}) \cup {\mathbf U}_{j'_{k+1}} \biggr)=\biggl({\mathbf U}_{j_1}\cap {\mathbf U}(j'_1,\cdots, j'_{k})\biggr)\bigcup \biggl({\mathbf U}_{j_1}\cap {\mathbf U}_{j'_{k+1}}\biggr).$$
Therefore there exists $g_1\in \Mor\bigl({\mathbf U}_{j_1}\cap {\mathbf U}(j'_1,\cdots, j'_{k})\bigr),$ and $g_2\in \Mor\bigl({\mathbf U}_{j_1}\cap {\mathbf U}_{j'_{k+1}}\bigr),$ such that,
$$f=g_2\circ g_1.$$
So, we have an isomorphism $g_0\in \Mor({\mathbf U}_{j'_k})\cap \Mor({\mathbf U}_{j'_{k+1}})$ satisfying,
\begin{equation}\label{E:divis}
\begin{split}
&g_1=g_0^{-1}\circ f'_{j'_k}\circ\cdots \circ f'_{j'_1},\\
&g_2=f'_{j'_{k+1}}\circ g_0.
\end{split}
\end{equation}
Since, $g_1\in \Mor\bigl({\mathbf U}_{j_1}\cap {\mathbf U}(j'_1,\cdots, j'_{k})\bigr),$ we can apply {\bf {[Ind assum 1]}} to obtain,
\begin{equation}\label{E:kto1}
\mbeta_{j_1}(g_1)=\mbeta_{j'_k}(g_0^{-1}\circ f'_{j'_k})\circ\cdots \circ \mbeta_{j'_1}(f'_{j'_1}).
\end{equation}
Similarly for $g_2,$ we have
\begin{equation}\label{E:kto2}
\mbeta_{j_1}(g_2)=\mbeta_{j'_{k+1}}(f'_{j'_{k+1}}\circ g_0).
\end{equation}
Composing \eqref{E:kto1} and \eqref{E:kto2} we arrive at our desired result,

\begin{equation}\label{E:finalk+1to1}
\begin{split}
\mbeta_{j_1}(f)&=\mbeta_{j_1}(g_2)\circ \mbeta_{j_1}(g_1)\\ 
&\hskip .2 cm [\hskip 0.05 cm \hbox{since}, {\mbeta}_{j_1}\hbox{ is a functor, and}\hskip 0.15 cm  g_2, g_1\in \Mor({\mathbf U}_{j_1})]\\
&=\mbeta_{j'_{k+1}}(f'_{j'_{k+1}}\circ g_0)\circ \mbeta_{j'_k}(g_0^{-1}\circ f'_{j'_k})\circ\cdots \circ \mbeta_{j'_1}(f'_{j'_1})\\
&=\mbeta_{j'_{k+1}}(f'_{j'_{k+1}})\circ \mbeta_{j'_k}(f'_{j'_k})\circ\cdots \circ \mbeta_{j'_1}(f'_{j'_1})\\
&[\hskip 0.05 cm \hbox{since},  {\mbeta}_{j'_{k+1}}|_{{\mathbf U}_{j'_{k+1} j'_{k}}}={\mbeta}_{j'_k}|_{{\mathbf U}_{j'_{k+1} j'_k}}].
\end{split}
\end{equation}
Hence proved.
\end{proof}
That means, if $f\in \Mor({\mathbf U}_{j_1}),$ and for $J'=\{j'_1,\cdots,j'_n\}\subset I,$ we have
\begin{equation}\nonumber
f=f'_{j'_n}\circ\cdots \circ f'_{j'_1},
\end{equation}
where $f'_{j'_k}\in \Mor({\mathbf U}_{j'_k}), k\in\{1, \cdots, n\},$ then, 
\begin{equation}\label{E:inpropglobalpsijij'}
\mbeta_{j_1}(f)=\mbeta_{j'_n}(f'_{j'_n})\circ\cdots\circ \mbeta_{j'_1}(f'_{j'_1}).
\end{equation}

Next we come to the second  induction assumption {\bf {[Ind assum 2]}}.
\begin{prop}\label{pr:glusheafk+1andn}
Suppose {\bf {[Ind assum 2]}} is true.  Let $\{j_1,\cdots, j_{k+1}\}\subset I, 
{J'}=\{j'_1,\cdots, j'_{n}\} \subset I,$ and 
$$f_{j_{k+1}}\circ\cdots\circ f_{j_1}=f'_{j'_{n}}\circ\cdots\circ f'_{j'_1}\in \Mor\biggl(\bigcup_{j_l\in J}{{\mathbf U}_{j_l}}\biggr)\bigcap \Mor\biggl(\bigcup_{j'_i\in J'}{\mathbf U}_{j'_i}\biggr),$$
where $f_{j_l}\in \Mor({\mathbf U}_{j_l}), l\in\{1,\cdots, k+1\}, f'_{j'_i}\in \Mor({\mathbf U}_{j'_i}), i\in\{1,\cdots, n\}.$ Then,
\begin{equation}\label{E:glusheafmandn}
\mbeta_{j_{k+1}}(f_{j_{k+1}})\circ\cdots\circ\mbeta_{j_1}(f_{j_1})=\mbeta_{j'_{n}}(f'_{j'_{n}})\circ \cdots\circ\mbeta_{j'_1}(f'_{j'_1}).
\end{equation}
\end{prop}
\begin{proof}
We follow our earlier notation, $${\mathbf U}_{j_1}\cup\cdots\cup{\mathbf U}_{j_k} \stackrel{\small{{\rm {notation}}}}{=}{\mathbf U}(j_1,\cdots, j_k),$$
and 
$${\mathbf U}_{j'_1}\cup\cdots\cup{\mathbf U}_{j'_n} \stackrel{\small{{\rm {notation}}}}{=}{\mathbf U}(j'_1,\cdots, j'_n).$$
Again  ${\mathbf U}({j_1},\cdots ,j_{k})\cap {\mathbf U}({j'_1,\cdots ,j'_n})$ and ${\mathbf U}_{j_{k+1}}\cap {\mathbf U}_{j'_{n}}$ are non empty [since, $s(f_{j_1})=s(f'_{j'_1}), t(f_{j_{k+1}})=t(f'_{j'_{n}}),$ and $s(f)\in \Obj({\mathbf U}_{j_1}), t(f)\in \Obj({\mathbf U}_{j_{k+1}})$].
Then we have,
\begin{equation}
\begin{split}
&f_{j_{k+1}}\circ \cdots \circ f_{j_1}=f'_{j'_n}\circ\cdots \circ f'_{j'_1}\\
&\in   \biggl({{\mathbf U}(j_1,\cdots, j_k)\cup{\mathbf U}_{j_{k+1}}}\biggr)\bigcap \biggl({\mathbf U}({j'_1,\cdots, j'_n}) \biggr)=\\
&\biggl({\mathbf U}_{j_{k+1}}\cap {\mathbf U}(j'_1,\cdots, j'_{n})\biggr)\bigcup \biggl({\mathbf U}({j_1,\cdots, j_k})\cap {\mathbf U}({j'_{1}},\cdots, j'_n \biggr).
\end{split}
\end{equation}
Therefore there exists $g_1\in \Mor\bigl({\mathbf U}({j_1\cdots j_k})\cap {\mathbf U}({j'_{1}}\cdots j'_n) \bigr),$ and $g_2\in \Mor\bigl({\mathbf U}_{j_{k+1}}\cap {\mathbf U}(j'_1,\cdots, j'_{n})\bigr),$ such that,
$$f'_{j'_n}\circ\cdots \circ f'_{j'_1}=f_{j_{k+1}}\circ \cdots \circ f_{j_1}=g_2\circ g_1.$$
So, we have a pair of  isomorphisms $g_0\in \Mor({\mathbf U}_{j_{k}})\cap \Mor({\mathbf U}_{j_{k+1}}), g'_0 \in \Mor({\mathbf U}_{j'_{n-1}})\cap \Mor({\mathbf U}_{j'_{n}}) $ satisfying,
\begin{equation}\nonumber
\begin{split}
&g_1=g_0^{-1}\circ f_{j_k}\circ\cdots \circ f_{j_1},\\
&g_2=f_{j_{k+1}}\circ g_0,
\end{split}
\end{equation}
and, 
\begin{equation}\nonumber
\begin{split}
&g_1={g'_0}^{-1}\circ f_{j'_{n-1}}\circ\cdots \circ f'_{j'_1},\\
&g_2=f_{j'_{n}}\circ g'_0.
\end{split}
\end{equation}
Now applying {\bf {[Ind assum 2]}}, and using same argument as previous proposition, we deduce ,
$$\mbeta_{j_{k+1}}(f_{j_{k+1}})\circ\cdots\circ\mbeta_{j_1}(f_{j_1})=\mbeta_{j'_{n}}(f'_{j'_{n}})\circ \cdots\circ\mbeta_{j'_1}(f'_{j'_1}).$$ Hence proved.
\end{proof}
By induction, we conclude, if $f_{j_l}\in \Mor({\mathbf U}_{j_l}), f'_{j'_i}\in \Mor({\mathbf U}_{j'_i})$ and for $ J=\{j_1, \cdots, j_m\}, J'=\{j'_1,\cdots,j'_n\}\subset I,$ we have
\begin{equation}\nonumber
f_{j_m}\circ\cdots \circ f_{j_1}=f'_{j'_n}\circ\cdots \circ f'_{j'_1},
\end{equation}
 then, 
\begin{equation}\label{E:inpropgl}
\mbeta_{j_m}(f_{j_m})\circ\cdots\circ\mbeta_{j_1}(f_{j_1})=\mbeta_{j'_n}(f'_{j'_n})\circ\cdots\circ \mbeta_{j'_1}(f'_{j'_1}).
\end{equation}
So,  \eqref{E:globalpsijij'} is proven.
\subsubsection*{{\bf{Proof of Theorem~\ref{Th:funcsheaf}}}} We conclude this subsection by presenting  the proof of  Theorem~\ref{Th:funcsheaf}. We have already derived all the required results. We only have to collect and organize them. We do not reiterate   the notations and conventions. They should be assumed to be as per with the Theorem~\ref{Th:funcsheaf} and subsequent part of this section.

Verification of {\bf {[Locality]}} condition  is virtually identical  to the given  in the proof of Proposition~\ref{pr:cover2sheaf}. For the sake of completeness, we only restate the result here. 

If, ${\mathbf \Psi}_1, {\mathbf \Psi}_2\in \Obj({\mathbf C}^{\mathbf U})={\rm Fun}(\mathbf U, \mathbf C),$ such that,
$${\mathbf \Psi}_1|_{{\mathbf U}_{\alpha}}= {\mathbf \Psi}_2|_{{\mathbf U}_{\alpha}}, \qquad\alpha\in I.$$
Then $\mbeta_1=\mbeta_2$ on $\mathbf U.$

Similarly, if $\bigl({\mathcal S}_1:\mbeta_1\Longrightarrow \mbeta_2\bigr), \bigl({\mathcal S}_2:\mbeta_1\Longrightarrow \mbeta_2\bigr)\in \Mor({\mathbf C}^{\mathbf U})={\mathcal N}(\mathbf U, \mathbf C),$ such that,
$${\mathcal S}_1|_{{\mathbf U}_{\alpha}}= {\mathcal S}_2|_{{\mathbf U}_{\alpha}}, \qquad\alpha\in I,$$
then ${\mathcal S}_1={\mathcal S}_2.$

Let us verify the {\bf {[Gluing]}} condition.

$\mbeta_{\alpha}\in \Obj({\mathbf C}^{{\mathbf U}_{\alpha}})={\rm Fun}({\mathbf U}_{\alpha}, {\mathbf C})$ given for each $\alpha\in I,$ such that:
$\mbeta_{\alpha}|_{{\mathbf U}_{\alpha \beta}}=\mbeta_{\beta}|_{{\mathbf U}_{\alpha \beta}},$ for every non-empty ${\mathbf U}_{\alpha \beta}.$
Then we define $(\mbeta:{\mathbf U}\longrightarrow {\mathbf C}) \in {\rm Fun}({\mathbf U}, {\mathbf C})=\Obj({\mathbf C}^{\mathbf U})$ as follows:
\begin{eqnarray}
&&\Obj({\mathbf U})\to \Obj({\mathbf C})\nonumber\\
&&\Obj(\bigcup_{\alpha\in I}{\mathbf U}_{\alpha})\to \Obj({\mathbf C})\nonumber\\
&& \hskip 1 cm a\mapsto \mbeta_{\alpha}(a), \qquad \hbox{if}\hskip 0.15 cm  a\in {\mathbf U}_{\alpha}, \alpha\in I,\label{E:gluobjgen}
\end{eqnarray}
and,
\begin{eqnarray}
&&\Mor({\mathbf U})\to \Mor({\mathbf C})\nonumber\\
&&{\mathbf {Gen}}\biggl(\bigcup_{\alpha \in I}\Mor({\mathbf U}_\alpha) \biggr)\to \Mor({\mathbf C}),\nonumber\\
&& f_{j_m}\circ\cdots\circ f_{j_1}\mapsto \mbeta_{j_m}(f_{j_m})\circ\cdots\circ \mbeta_1(f_{j_1}), \nonumber\\
&&\hbox{where}\hskip 0.15 cm f_{j_k}\in \Mor({\mathbf U}_{j_k}), k\in \{1,\cdots , m\},\nonumber \\
&&\hbox{and,}\hskip 0.15 cm J:=\{j_1,\cdots j_m\}\subset I.\label{E:glumorgen}
\end{eqnarray}
Since $\mbeta_{\alpha}|_{{\mathbf U}_{\alpha \beta}}=\mbeta_{\beta}|_{{\mathbf U}_{\alpha \beta}},$ \eqref{E:gluobjgen} is well defined. Also  \eqref{E:globalpsijij'} ensures that \eqref{E:glumorgen} is well defined. Thus we have a 
$(\mbeta:{\mathbf U}\longrightarrow {\mathbf C}) \in {\rm Fun}({\mathbf U}, {\mathbf C})$ satisfying,
$$\mbeta_{\alpha}=\mbeta|_{\alpha}, \alpha\in I.$$

Next, suppose $\bigl({\mathcal S}_{\alpha}:{\mbeta}_{\alpha}\Longrightarrow \mbeta'_{\alpha}\bigr)\in \Mor({\mathbf C}^{{\mathbf U}_{\alpha}})={\mathcal N}({\mathbf U}_{\alpha}, {\mathbf C})$ is given for each $\alpha\in I,$ such that,
$\{\mbeta_{\alpha}\}$ glue to form a $\mbeta\in {\rm Fun}(\mathbf U, \mathbf C),$ and  similarly, $\{\mbeta'_{\alpha}\}$ glue to form a $\mbeta'\in {\rm Fun}(\mathbf U, \mathbf C),$ [as described in the previous part of the proof] and ${\mathcal S}_{\alpha}$ satisfy 
\begin{equation}\label{E:natovgen}
{\mathcal S}_{\alpha}|_{{\mathbf U}_{\alpha \beta}}={\mathcal S}_{\beta}|_{{\mathbf U}_{\alpha \beta}}.
\end{equation}
We define ${\mathcal S}\in {\mathcal N}({\mathbf U}, {\mathbf C})$ by:
\begin{equation}\label{E:glonatgen}
{\mathcal S}(a)={\mathcal S}_{\alpha}(a), \hskip 0.4 cm \hbox{for}\hskip 0.15 cm a\in \Obj({\mathbf U}_{\alpha}), \alpha\in I.
\end{equation}
Above equation is well defined because of \eqref{E:natovgen}. The proof that \eqref{E:glonatgen} defines a natural transformation from $\mbeta$ to $\mbeta',$ is exactly same as the counter part in Proposition~\ref{pr:cover2sheaf}.

And, that completes the proof of Theorem~\ref{Th:funcsheaf}.$\fbox{}$

\section{Sheaves of categorical groups}\label{S:sheavcaygrop} 
Objective of this section  is to define the (pre)sheaves of categorical groups, and construct an example for the same. In section~\ref{S:mbsieve} we made some passing remarks about categorical group valued presheaves. Here we will take a more formal approach. Before that let us briefly review categorical groups, and associated notions.

\subsection{Categorical groups}
There are many equivalent definitions of a \textit{categorical group}  available in literature [17, 20, 24, 26]. For our purpose we will mostly think of a categorical group in terms of a \textit{crossed-module}. A {\em categorical group} $\mbg$ is given by a category  ${\mathcal G}$  along with a functor
\begin{equation}\label{E:compcatG}
\mbg\times \mbg\to\mbg 
\end{equation}
which makes both $\Obj(\mbg)$ and $\Mor(\mbg)$ groups. Some of the immediate consequences are as follows: 
\begin{itemize}
	\item[]{(i) the identity-assigning map 
$$\Obj(\mbg)\to\Mor(\mbg): x\mapsto 1_x$$
is a homomorphism, and so, $1_e$ is the identity element in $\Mor(\mbg),$ where $e$ is the identity element in $\Obj(\mbg),$   } 
\item[]{(ii) the source and target maps
\begin{equation}\label{E:st}
s,t:\Mor(\mbg)\to \Obj(\mbg) 
\end{equation}
are both homomorphisms;}
\item[]{(iii) functoriality of (\ref{E:compcatG}) implies the following exchange law, 
\begin{equation}\label{E:exch}
(\phi_2\psi_2)\circ(\phi_1\psi_1)=(\phi_2\circ \phi_1)(\psi_2\circ \psi_1),
\end{equation}
whenever right hand side is well defined for $\phi_1, \phi_2, \psi_1, \psi_2\in \Mor(\mbg)$. Here and onwards $\circ$ will denote the composition of morphisms (as usual) and juxtaposition of two elements denote group product.} 
\end{itemize}
Specializing to the case, where  both $\Obj(\mbg)$ and $\Mor(\mbg)$ are Lie groups and the maps $s$, $t$ and $x\mapsto 1_x$ are smooth, we have a {\em categorical Lie group} $\mbg$.

A {\em crossed-module} is given by a pair of groups $G$ and $H$, along with maps
$$\alpha:G\times H\to H:(g, h)\mapsto \alpha_g(h)\quad\hbox{and}\quad \tau:H\to G,$$
where $\tau$ is a homomorphism,  $\alpha_g$ is an automorphism of $H$ for each $g\in G$, and the map $g\mapsto \alpha_g\in {\rm Aut}(H)$ is a homomorphism.  The  map $\tau$ and the map $\alpha$ interrelated via following identities
\begin{equation}\label{liecross0}
	\begin{split}
		&\tau(\alpha_g(h))=g\tau(h)g^{-1},\\
  &\alpha_{\tau(h)}(h')=hh'h^{-1}\qquad\hbox{for all $g\in G, h, h'\in H$.}
\end{split}
\end{equation}
We write a crossed module as $(G, H, \alpha,\tau).$ When $G$ and $H$ are Lie groups, and $\alpha$ and $\tau$ are smooth,   $(G, H, \alpha,\tau)$ is called a {\em Lie crossed module}. For us, most useful property of a (Lie) crossed module is the following.

It is well known that there is a one-to-one correspondence between categorical (Lie) groups and (Lie) crossed modules [5, 6, 16]. The bijection is given as follows. 

Let $\mbg$ be a categorical group. We take $G:=\Obj(\mbg)$, $H:=\ker s$, $\tau=t|_H$, and
$$\alpha_g(h)=1_gh1_g^{-1}$$ for all $g\in G$ and $h\in H$. Then we have a group isomorphism $\Mor(\mbg)\stackrel{\simeq}\to H\rtimes_{\alpha}G$ defined by the map
\begin{equation}\nonumber
	\begin{split}
		\Mor(\mbg)&\to H\rtimes_{\alpha}G\\
		&:\phi\mapsto \bigl(\phi1_{s(\phi)^{-1}}, s(\phi) \bigr).
 \end{split}
\end{equation}
The target map $t$, viewed as a mapping $H\rtimes_{\alpha}G\to G$, is given by
\begin{equation}\label{E:ttau}
t(h,g)=\tau(h)g\qquad\hbox{for all $(h,g)\in H\rtimes_{\alpha}G$.}
\end{equation}
We note here that the group product in $H\rtimes_{\alpha}G$ is given by the usual group multiplication law for a semi direct product
\begin{equation}\label{E:hgprod}
(h_2,g_2)(h_1,g_1)=\bigl(h_2\alpha_{g_2}(h_1), g_2g_1\bigr)
\end{equation}
for all $(h_2,g_2), (h_1,g_1)\in H\rtimes_{\alpha}G$. It also easily follows that the composition of morphisms in $\Mor(\mbg)\simeq H\rtimes_{\alpha}G$ is given by
\begin{equation}\label{E:hgmorcomp}
(h_2,g_2)\circ (h_1, g_1)=(h_2h_1,g_1);
\end{equation}
of course this composition is defined only when 
$$\tau(h_1)g_1=t(h_1,g_1)=s(h_2,g_2)=g_2.$$
A \textit{morphism between a pair of categorical groups} ${\mathcal G}$ and ${\mathcal H}$ is a functor [5, 6]
\begin{equation}\label{E:homcatgroups}
\lambda: {\mathcal G}\longrightarrow {\mathcal H},
\end{equation}
such that, both $\lambda:\Obj(\mathcal G)\to \Obj({\mathcal H})$  and $\lambda:\Mor(\mathcal G)\to \Mor({\mathcal H})$ are homomorphisms of groups.

The category of all categorical groups, that is, the category whose objects are categorical groups and morphisms are morphisms between categorical groups (defined in \eqref{E:homcatgroups}), will be denoted as
$${\mathbf {CatGrp}}.$$
Clearly  there is a full, faithful, forgetful functor 
$$\mathbf {CatGrp}\longrightarrow \mathbf {Cat},$$
and ${\mathbf {CatGrp}}$ is a full subcategory of $\mathbf {Cat}.$ Our next goal is to consider the ${\mathbf {CatGrp}}\subset {\mathbf {Cat}}$ valued presheaf over a category ${\mathcal C}$ of a collection of small categories. 
We define a presheaf of {categorical groups} or a $\mathbf {CatGrp}$-valued presheaf, over $\mathcal C,$ to be a contravariant functor from $\mathcal C$ to $\mathbf {CatGrp}:$
\begin{equation}\label{E:catgrpprshf}
{\rho}: {\mathcal C}^{\rm op}\to \mathbf {CatGrp}.
\end{equation}\label{E:catofcatgrpprschf}
We will denote the \textit{category of presheaves of categorical groups} by \begin{equation}
{{\mathbf {Prsh}}}(\mathcal C, {\mathbf {CatGrp}}).
\end{equation}

Let $\mathbf B$ be a topological category and ${\widetilde {\mathcal O}}(\mathbf B)$ be as defined in \eqref{E:genobopen}. A presheaf of categorical groups over ${\widetilde {\mathcal O}}(\mathbf B)$  is a \textit{sheaf of categorical groups} over ${\widetilde {\mathcal O}}(\mathbf B)$     if it satisfies {\bf {[Locality]}} and {\bf {[Gluing]}} conditions listed between \eqref{E:objsheafcov}--\eqref{E:ovmor}.

\subsection{Functor category $\mbg^\mathbf U$}\label{s:fc}
Let $\mathbf U$ be a category, and $\mbps_1, \mbps_2:\mathbf U\to \mbg$ be a pair of functors from $\mathbf U$ to a categorical group $\mbg.$ Then the pointwise product $\mbps_2\mbps_1:\mathbf U\to \mbg$ is also a functor. In other words, we have following result [Proposition 3.2, [15]].
\begin{lemma}\label{l:grpfunc}
The set of all functors from $\mathbf U$ to $\mbg$ form a group.	
\end{lemma}
\begin{proof}
	$\mbps_2\mbps_1$ is defined as follows. \begin{equation}\label{E:prodfun}
	\bigl(\mbps_2\mbps_1\bigr)(x)=\mbps_2(x)\mbps_1(x),\end{equation}
	where $x$ is in $\Obj(\mathbf U)$ or $\Mor(\mathbf U).$ Thus, if $f_2, f_1\in \Mor(\mathbf U)$ and they are composable, then
	\begin{eqnarray}
		(\mbps_2\mbps_1)(f_2\circ f_1)&&= \mbps_2(f_2\circ f_1)\mbps_1(f_2\circ f_1)\nonumber\\
		 &&=\bigl(\mbps_2(f_2)\circ \mbps_2(f_1)\bigr)\bigl(\mbps_1(f_2)\circ \mbps_2(f_1)\bigr)\nonumber\\
   &&=\bigl(\mbps_2(f_2) \mbps_1(f_2)\bigr)\circ\bigl(\mbps_2(f_1)\mbps_1(f_1)\bigr)\nonumber[\hbox{using \eqref{E:exch}}]\\
                 &&=\bigl(\mbps_2\mbps_1\bigr)(f_2)\circ\bigl(\mbps_2\mbps_1\bigr)(f_1).\nonumber
 	\end{eqnarray}
	So $\mbps_2\mbps_1$ is a functor. It is obvious that the constant functor 
	\begin{eqnarray}
		&\mbps_0:&\mathbf U\to \mbg\nonumber\\
			      &&a\mapsto e,\nonumber\\
		&&f\mapsto 1_e\nonumber
	\end{eqnarray}
	defines the multiplicative identity, where $a\in \Obj(\mathbf U)$ { and} $f\in \Mor(\mathbf U),$ and $e, 1_e$ respectively denote the group identity elements in $\Obj(\mathbf U)${ and} $\Mor(\mathbf U).$ Group inverse $\mbps^{\rm inv}$ is given by $\mbps^{\rm inv}(x)=\mbps(x)^{-1},$ where $x\in \Obj(\mbg)$ or $\Mor(\mbg).$
\end{proof}
Now suppose $\mathcal T:\mbps_1\Longrightarrow \mbps_2$ and ${\mathcal T}':\mbps'_1\Longrightarrow \mbps'_2$ are a pair of natural transformations between respective functors in $\{\mathbf U\to \mbg\};$ that is, the diagrams in \eqref{E:Tnatural2} commute, where $a\stackrel{f}{\rightarrow} b$ is a morphism in category $\mathbf U.$
\begin{equation} \label{E:Tnatural2}
 \xymatrixcolsep{5pc}\xymatrix{
         \ar[d]^{{\mathcal T}(a)} \mbps_1(a)       \ar[r]^-{\mbps_1(f)} & \mbps_1(b) \ar[d]_{{\mathcal T}(b)} \\
\mbps_2(a) \ar[r]_-{\mbps_2(f)}& \mbps_2(b) 
}
\hskip 1in
\xymatrix{
         \ar[d]^{{\mathcal T}'(a)} \mbps'_1(a)       \ar[r]^-{\mbps'_1(f)} & \mbps'_1(b) \ar[d]_{{\mathcal T}'(b)} \\
\mbps'_2(a) \ar[r]_-{\mbps'_2(f)}& \mbps'_2(b) 
}
\end{equation}
We can define a product natural transformation, given by \begin{equation}\label{E:prodnat}
{\mathcal T}'\mathcal T(a):={\mathcal T}'(a){\mathcal T}(a),\end{equation} for any $a\in \Obj(\mathbf U).$ Functoriality of the group product $\mbg\times \mbg\to \mbg$ ensures that ${\mathcal T}'\mathcal T$ is a well defined natural transformation; that is, following diagram commutes:
\begin{equation} \label{E:Tnatural3}
\xymatrix{
         \ar[d]^{{\mathcal T}'(a){\mathcal T}(a)} \mbps_1(a)  \mbps'_1(a)     \ar[rr]^-{\mbps_1(f)\mbps'_1(f)} & &\mbps_1(b)\mbps'_1(b) \ar[d]_{{\mathcal T}'(b){\mathcal T}(b)} \\
\mbps_2(a)\mbps'_2(a) \ar[rr]_-{\mbps_2(f)\mbps'_2(f)}& & \mbps_2(b) \mbps'_2(b) 
}.
\end{equation}
For any ${\mathcal T}$ we have a corresponding (multiplicative) inverse given by
$${\mathcal T}^{\rm inv}(a):=({\mathcal T}(a))^{-1}.$$

Moreover, if ${\mathcal T}_1:\mbps_1\Longrightarrow \mbps_2$ and ${\mathcal T}_2:\mbps_2\Longrightarrow \mbps_3$ are natural transformations, then we have a composite natural transformation ${\mathcal T}_2\circ{\mathcal T}_1:\mbps_1\Longrightarrow \mbps_3$,
$$({\mathcal T}_2\circ{\mathcal T}_1)(a)={\mathcal T}_2(a)\circ{\mathcal T}_1(a), \qquad [\hbox{for all}\hskip 0.2 cm a\in \Obj(\mathbf U)].$$ 
Also, using functoriality of $\mbg\times \mbg\to \mbg$, it is easy to show that, 
\begin{equation}\label{E:UGopfunct}
({\mathcal T}'_2{\mathcal T}_2)\circ ({\mathcal T}'_1{\mathcal T}_1)= ({\mathcal T}'_2\circ {\mathcal T}'_1)({\mathcal T}_2\circ {\mathcal T}_1),
\end{equation}
when ${\mathcal T}_2\circ {\mathcal T}_1$ and ${\mathcal T}'_2\circ {\mathcal T}'_1$ are defined. The natural transformation ${\mathcal T}_0:\mbps_0\Longrightarrow \mbps_0$ defined as $${\mathcal T}_0(a)=1_{e}, \qquad\forall a\in \Obj(\mathbf U)$$
is the multiplicative identity.

Thus we make following proposition [Proposition 3.4, [15]]:
\begin{prop}\label{P:funccat}
	Let $\mbg^{\mathbf U}$ be the  category of functors from $\mathbf U$ to $\mbg$:
	\begin{equation}\label{E:funcat}
		\begin{split}
		&\Obj(\mbg^{\mathbf U}):={\rm Fun}(\mathbf U, \mbg)\\
		       &\Mor(\mbg^{\mathbf U}):={\mathcal N}(\mathbf U, \mbg).
		\end{split}
	\end{equation}
	Then $\mbg^{\mathbf U}$ is a categorical group.
\end{prop}

Now suppose,  as in Theorem~\ref{Th:funcsheaf},  $\widetilde {\mathcal O}(\mathbf B)$ is the category of open subcategories of the topological groupoid $\mathbf B,$    and    $\{{\mathbf U}_{\alpha}\}$  is an open categorical cover of ${\mathbf U}.$ Then by Lemma~\ref{L:escobjmor}, for any functor  $\biggl({\mathbf U}\xrightarrow{\mathbf \Theta}{\mathbf V}\biggr)\in \Mor\biggl({\widetilde {\mathcal O}}(\mathbf B)\biggr),$ we have a functor
${\mathbf \rho}(\mbl) : {\mathcal G}^{\mathbf V}\longrightarrow {\mathcal G}^{\mathbf U}:$
\begin{equation}\label{E:secobjmorgrp}
\begin{split}
{\mathbf \rho}(\mbl) :  &\Obj({\mathcal G}^{\mathbf V})\longrightarrow \Obj({\mathcal G}^{\mathbf U}),\\
&{\mathbf \Phi}\mapsto {\mathbf \Phi}\mbl,\\
{\mathbf \rho}(\mbl) :  &\Mor({\mathcal G}^{\mathbf V})\longrightarrow \Mor({\mathcal G}^{\mathbf U}).\\
&{\mathcal T}\mapsto {\mathcal T}\mbl.
\end{split}
\end{equation}
In fact, the functor in \eqref{E:secobjmorgrp} is a morphism of categorical groups defined in \eqref{E:homcatgroups}; that is, the functor defines a pair of group homomorphisms,
\begin{equation} \label{E:zz}
\begin{split} 
&\Obj\bigl({\mathcal G}^{\mathbf V}\bigr)\to \Obj\bigl({\mathcal G}^{\mathbf U}\bigr), \hskip 0.1 cm \hbox{and}\\
&\Mor\bigl({\mathcal G}^{\mathbf V}\bigr)\to \Mor\bigl({\mathcal G}^{\mathbf U}\bigr).
\end{split}
\end{equation}
\begin{lemma}\label{L:catgrpprschf}
Let ${\mathbf U}, {\mathbf V}\subset {\mathbf B},$ and $\mathcal G$ be a fixed categorical group. Let ${\mathcal G}^{\mathbf U}, {\mathcal G}^{\mathbf V}$  be the respective categorical groups of functors defined in Proposition~\ref{P:funccat}. Then for any $\biggl({\mathbf U}\xrightarrow{\mathbf \Theta}{\mathbf V}\biggr)\in \Mor\biggl({\widetilde {\mathcal O}}(\mathbf B)\biggr),$ the functor ${\mathbf \rho}(\mbl)$, defined in \eqref{E:secobjmorgrp}, is a morphism from ${\mathcal G}^{\mathbf V}$ to ${\mathcal G}^{\mathbf U}.$
\end{lemma}
\begin{proof}
Proof directly follows from the group product  defined for objects and  morphisms respectively in \eqref{E:prodfun} and\eqref{E:prodnat}. Let ${\mathbf \Phi}, {\mathbf \Phi}'\in \Obj({\mathcal G}^{\mathbf V}).$ Then for any $x\in \Obj(\mathbf U) \hskip 0.1 cm \hbox{or} \hskip 0.1 cm \Mor(\mathbf U),$ we have
\begin{equation}
\begin{split}
&\biggl(\rho (\mbl)({\mathbf \Phi}{\mathbf \Phi}')\biggr)(x)=\bigl({\mathbf \Phi}{\mathbf \Phi}'\mbl\bigr)(x)\hskip 0.2 cm \hbox{[using \eqref{E:secobjmorgrp}]}\\
&=\bigl({\mathbf \Phi}{\mathbf \Phi}'\bigr)(\mbl(x))\\
&= {\mathbf \Phi}(\mbl(x)){\mathbf \Phi}'(\mbl(x)) \hskip 0.2 cm \hbox{[using \eqref{E:prodfun}]}\\
&=\biggl(\rho(\mbl)({\mathbf \Phi})\biggr)(x)\biggl(\rho(\mbl)({\mathbf \Phi}')\biggr)(x)\hskip 0.2 cm \hbox{[using \eqref{E:secobjmorgrp}]}\\
&= \biggl(\rho(\mbl)({\mathbf \Phi})\rho(\mbl)({\mathbf \Phi}')\biggr)(x) \hskip 0.2 cm \hbox{[using \eqref{E:prodfun}]}.
\end{split}
\end{equation}
So,$$\rho(\mbl)({\mathbf \Phi})\rho(\mbl)({\mathbf \Phi}')=\rho (\mbl)({\mathbf \Phi}{\mathbf \Phi}').$$
Similarly, if ${\mathcal T}, {\mathcal T}'\in \Mor({\mathcal G}^{\mathbf V}),$ the using \eqref{E:prodnat} we can show
$$\rho(\mbl)({\mathcal T})\rho(\mbl)({\mathcal T}')=\rho (\mbl)({\mathcal T}{\mathcal T}').$$
Thus $\rho(\mbl)$ defines a pair of group homomorphisms
\begin{equation} \nonumber
\begin{split} 
&\Obj\bigl({\mathcal G}^{\mathbf V}\bigr)\to \Obj\bigl({\mathcal G}^{\mathbf U}\bigr), \hskip 0.1 cm \hbox{and}\\
&\Mor\bigl({\mathcal G}^{\mathbf V}\bigr)\to \Mor\bigl({\mathcal G}^{\mathbf U}\bigr).
\end{split}
\end{equation}
And, the functor $\rho(\mbl),$ as per definition in \eqref{E:homcatgroups}, is a morphism between categorical groups.

\end{proof}

In other words
\begin{equation}\label{E:presheafsecgrp}
\begin{split}
{\rho} : &\bigl({\widetilde {\mathcal O}}({\mathbf B})\bigr)^{\rm op}\to {\mathbf {CatGrp}},\\
&\Obj\biggl({\widetilde {\mathcal O}}({\mathbf B})\biggr)\to \Obj({\mathbf {CatGrp}}),\\
&{\mathbf U}\mapsto {\mathcal G}^{\mathbf U},\\
&\Mor\biggl({\widetilde {\mathcal O}}({\mathbf B})\biggr)\to \Mor(\mathbf {CatGrp}),\\
&\biggl({\mathbf U}\xrightarrow{\mbl}{\mathbf V}\biggr)\mapsto \biggl({\mathcal G}^{\mathbf V}\xrightarrow{{{\rho}(\mbl)}} {\mathcal G}^{\mathbf U}\biggr),
\end{split}
\end{equation}
is an element of ${{\mathbf {Prsh}}}(\mathcal C, {\mathbf {CatGrp}}).$

Theorem~\ref{Th:funcsheaf} implies that  $\mathbf {CatGrp}$-valued presheaf $\mathbf \rho$ in \eqref{E:presheafsecgrp} is actually a  $\mathbf {CatGrp}$-valued sheaf. Hence we conclude,
\begin{prop}\label{Pr:funcsheafgrp}
Let $\widetilde {\mathcal O}(\mathbf B)$ be the category of open subcategories for a  topological groupoid $\mathbf B.$ Let $\mathbf \rho$ be as defined in \eqref{E:presheafsecgrp}. Then $\mathbf \rho$ is a $\mathbf {CatGrp}$-valued sheaf over $\widetilde {\mathcal O}(\mathbf B).$
\end{prop}
\section*{appendix}

\subsection*{Proof of Proposition~\ref{pr:boolean}}

Let ${\mathbf U}, {\mathbf V}, {\mathbf W}$ be subcategories of a category $\mathbf B.$ Let  intersection and union  of a pair of subcategories be respectively as defined  in \eqref{E:objmorint}, \eqref{E:objmoruni}. If $\mathbf U,$ $\mathbf V$ and $\mathbf W$ are mutually disjoint, then the statements in Proposition~\ref{pr:boolean} are same as standard set theoretic statements, because
$$\Mor(\mathbf U\cup \mathbf V)=\Mor(\mathbf U)\cup \Mor(\mathbf V), \qquad \hbox{if}\hskip 0.1 cm \mathbf U\cap \mathbf V=\emptyset.$$
If one of the $\mathbf U, \mathbf V, \mathbf W$ is disjoint with other two, then also the proof is straightforward. In what follows, we will assume
\begin{equation}\nonumber
\begin{split}
&\mathbf U\cap \mathbf V\neq \emptyset,\\
&\mathbf U\cap \mathbf W\neq \emptyset,\\
&\mathbf W\cap \mathbf V\neq \emptyset.
\end{split}
\end{equation}

\subsubsection*{\bf{Proof of identity (i)}}
The second equation of (i),  
$${\mathbf U}\cap({\mathbf V} \cap {\mathbf W})=({\mathbf U}\cap {\mathbf V}) \cap {\mathbf W},$$
is an obvious consequence of the definition of intersection in \eqref{E:objmorint}. Let us prove the first equation of (i). 

We first  classify  morphisms in ${\mathbf U}\cup ({\mathbf V}\cup {\mathbf W}).$ If $f\in \Mor\biggl({\mathbf U}\cup ({\mathbf V}\cup {\mathbf W})\biggr),$ then  \eqref{E:objmoruni} implies following, not necessarily mutually exclusive, classification (see left hand side of Figure 3) :
\begin{figure}[h]
\begin{center}
\includegraphics[height=3.27in,width=6 in]{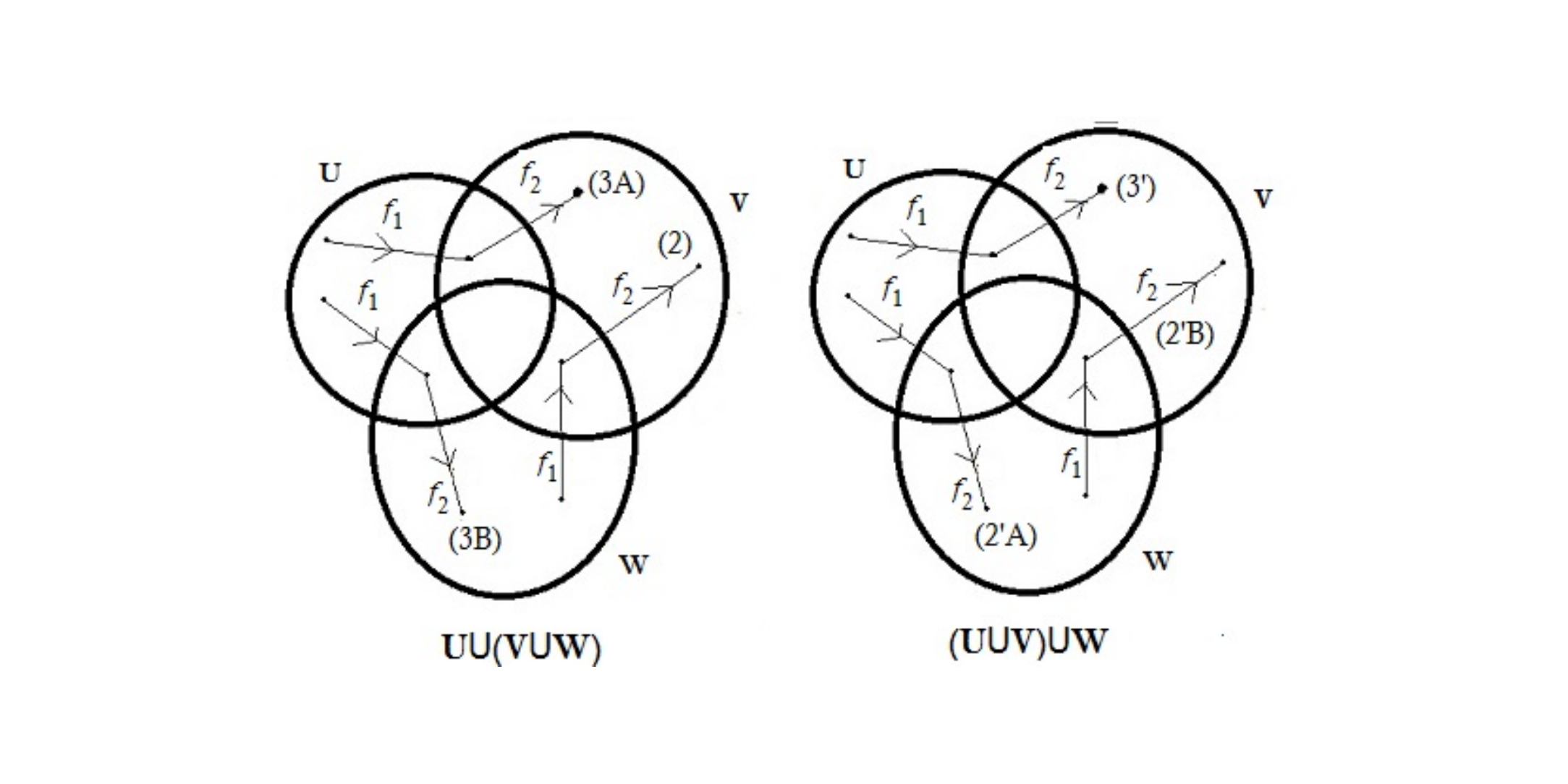}
\caption{ Comparison between morphisms in ${\mathbf U}\cup ({\mathbf V}\cup {\mathbf W})$ and $({\mathbf U}\cup {\mathbf V})\cup {\mathbf W}$ }
\end{center}
\end{figure}\label{f:appeni}

\begin{itemize}
\item[]{(1) $f\in \Mor({\mathbf U})\cup \Mor({\mathbf V})\cup \Mor({\mathbf W})$.
}
\item[]{(2)  $f$ is of the form $f_2\circ f_1,$ where $f_2, f_1\in \Mor({\mathbf V})\cup\Mor({\mathbf W}),$ and  $s(f_2)=t(f_1)\in \Obj({\mathbf V}\cap {\mathbf W})=\Obj({\mathbf V})\cap \Obj({\mathbf W}).$
}
\item[]{(3) $f$ is of the form $f_2\circ f_1,$ where $f_2, f_1\in \Mor({\mathbf U})\cup\Mor({\mathbf V}\cup{\mathbf W}),$ and $s(f_2)=t(f_1)\in \Obj\bigl({\mathbf U}\cap ({\mathbf V}\cup{\mathbf W})\bigr).$ By \eqref{E:objmorint}, \eqref{E:objmoruni}, for the objects, we have following identity
$\Obj\bigl({\mathbf U}\cap ({\mathbf V}\cup{\mathbf W})\bigr)= \Obj\bigl({\mathbf U}\bigr)\cap \biggl(\Obj\bigl({\mathbf V}\bigr)\cup\Obj\bigl({\mathbf W}\bigr)\biggr)
=\biggl(\Obj\bigl({\mathbf U}\bigr)\cap\Obj\bigl({\mathbf V}\bigr)\biggr)\cup \biggl(\Obj\bigl({\mathbf U}\bigr)\cap\Obj\bigl({\mathbf W}\bigr)\biggr).$ So, we may further classify the morphisms in (3) as,
\begin{itemize}
\item[]{(3A) $f=f_2\circ f_1$ such that $f_2, f_1\in \Mor(\mathbf U)\cup \Mor(\mathbf V)$ and $s(f_2)=t(f_1)\in \biggl(\Obj\bigl({\mathbf U}\bigr)\cap\Obj\bigl({\mathbf V}\bigr)\biggr).$}
\item[]{(3B) $f=f_2\circ f_1$ such that $f_2, f_1\in \Mor(\mathbf U)\cup\Mor(\mathbf W)$ and $s(f_2)=t(f_1)\in \biggl(\Obj\bigl({\mathbf U}\bigr)\cap\Obj\bigl({\mathbf W}\bigr)\biggr).$   }
\end{itemize}
}
\end{itemize}

\begin{figure}[h]
\begin{center}
\includegraphics[height=2.5 in,width=5.3 in]{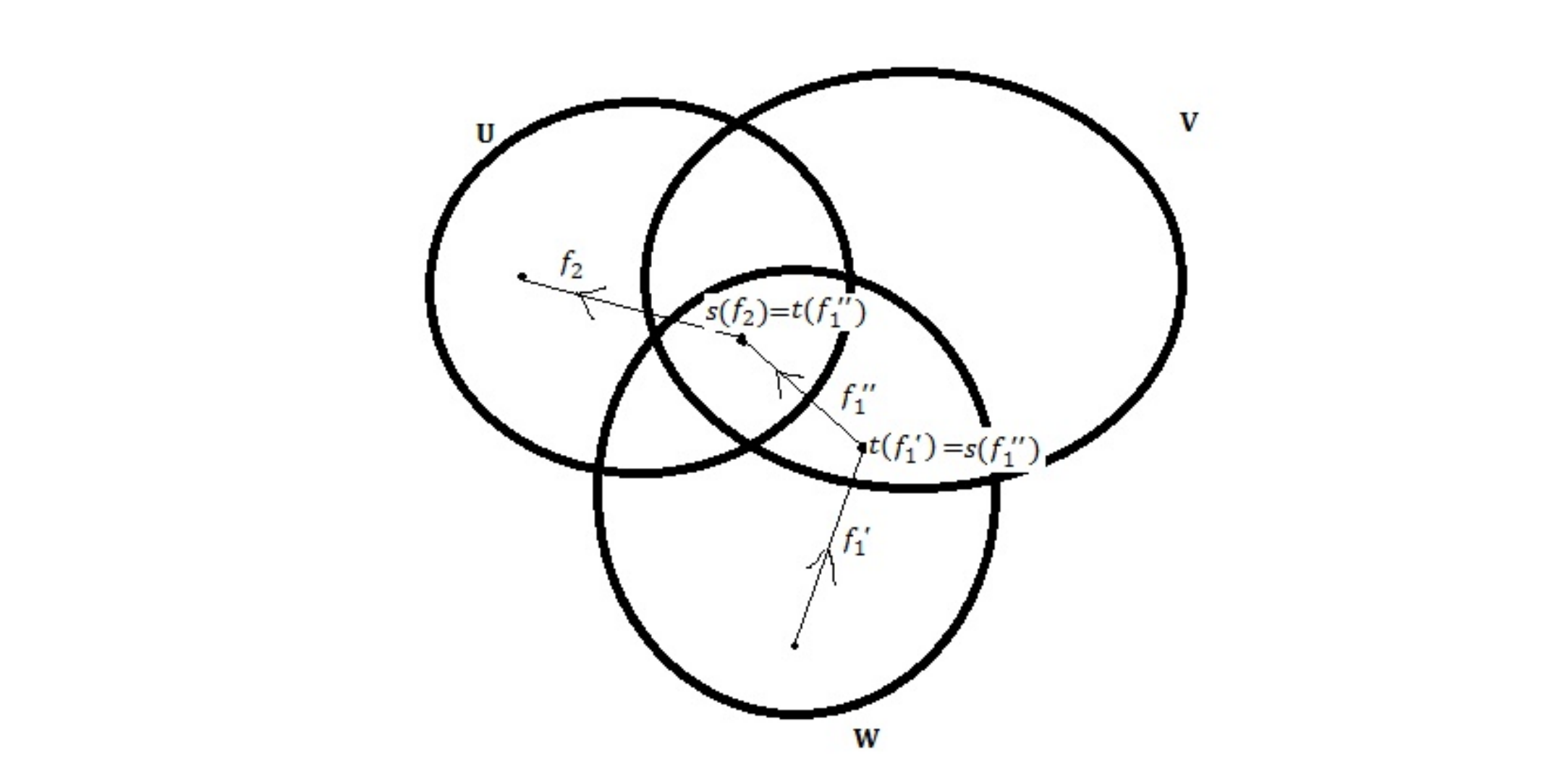}
\caption{ }
\end{center}
\end{figure}\label{f:ommited}

[It may appear that in (3) we have overlooked the case when $f''_1, f'_1\in \Mor(\mathbf V)\cup  \Mor(\mathbf W)$ and $f_1=f''_1\circ f'_1, s(f''_1)=t(f'_1)\in \Obj({\mathbf V}\cap {\mathbf W}).$ But, note that in that case $t(f''_1)=s(f_2)\in \Obj\biggl({\mathbf U}\cap {\mathbf V}\biggr)\cap\Obj\biggl({\mathbf V}\cap {\mathbf W}\biggr)=\Obj\biggl({\mathbf U}\cap {\mathbf W}\biggr)\cap\Obj\biggl({\mathbf V}\cap {\mathbf W}\biggr)=\Obj({\mathbf U}\cap{\mathbf V}\cap {\mathbf W}),$ and therefore $f''_1\circ f'_1=f_1\in \Mor(\mathbf W)$ or $f''_1\circ f'_1=f_1\in \Mor(\mathbf V).$ Such  scenarios have already been taken care of respectively by (3B) and (3A). Figure 4 illustrates the case when $f''_1\circ f'_1=f_1\in \Mor(\mathbf W).$]

We proceed with the  classification of  morphisms in $({\mathbf U}\cup {\mathbf V})\cup {\mathbf W}.$ If $f\in \Mor\biggl(\bigl({\mathbf U}\cup {\mathbf V}\bigr)\cup {\mathbf W}\biggr),$ then  we obtain following classification  (see right hand side of  Figure 3) :
\begin{itemize}
\item[]{(\'1) $f\in \Mor({\mathbf U})\cup \Mor({\mathbf V})\cup \Mor({\mathbf W})$.
}
\item[]{(\'2) $f$ is of the form $f_2\circ f_1,$ where $f_2, f_1\in \Mor\biggl(\bigl({\mathbf U}\cup {\mathbf V}\bigr)\cap{\mathbf W}\biggr)$ and $s(f_2)=t(f_1)\in \Obj\biggl(\bigl({\mathbf U}\cup {\mathbf V}\bigr)\cap{\mathbf W}\biggr).$ Using the identity
$ \Obj\biggl(\bigl({\mathbf U}\cup {\mathbf V}\bigr)\cap{\mathbf W}\biggr)
=\biggl(\Obj\bigl({\mathbf U}\bigr)\cap\Obj\bigl({\mathbf W}\bigr)\biggr)\cup \biggl(\Obj\bigl({\mathbf V}\bigr)\cap\Obj\bigl({\mathbf W}\bigr)\biggr),$ we further classify  morphisms in (\' 2) as,
\begin{itemize}
\item[]{(\'2A) $f=f_2\circ f_1$ such that $f_2, f_1\in \Mor(\mathbf U)\cup \Mor(\mathbf W)$ and $s(f_2)=t(f_1)\in \biggl(\Obj\bigl({\mathbf U}\bigr)\cap\Obj\bigl({\mathbf W}\bigr)\biggr).$}
\item[]{(\'2B) $f=f_2\circ f_1$ such that $f_2, f_1\in \Mor(\mathbf V)\cup \Mor(\mathbf W)$ and $s(f_2)=t(f_1)\in \biggl(\Obj\bigl({\mathbf V}\bigr)\cap\Obj\bigl({\mathbf W}\bigr)\biggr).$   }
\end{itemize}
}

\item[]{(\' 3) $f=f_2\circ f_1,$ where $f_2, f_1\in \Mor(\mathbf U)\cup \Mor(\mathbf V)$ and $s(f_2)=t(f_1)\in \biggl(\Obj\bigl({\mathbf U}\bigr)\cap\Obj\bigl({\mathbf V}\bigr)\biggr).$
}
\end{itemize}
We make following correspondence between two classifications:
\begin{equation}\nonumber
\begin{split}
&\rm{(1)}\Longleftrightarrow \hbox{(\'1)}\\
&\rm{(2)}\Longleftrightarrow \hbox{(\'2B)}\\
&\rm{(3A)}\Longleftrightarrow \hbox{(\'3)}\\
&\rm{(3B)}\Longleftrightarrow \hbox{(\'2A)}.
\end{split}
\end{equation}
Hence we conclude
$${\mathbf U}\cup ({\mathbf V}\cup {\mathbf W})=({\mathbf U}\cup {\mathbf V})\cup {\mathbf W}.$$
\subsection*{Proof of identity (ii)}

\begin{figure}[h]
\begin{center}
\includegraphics[height=2.5in,width=5in]{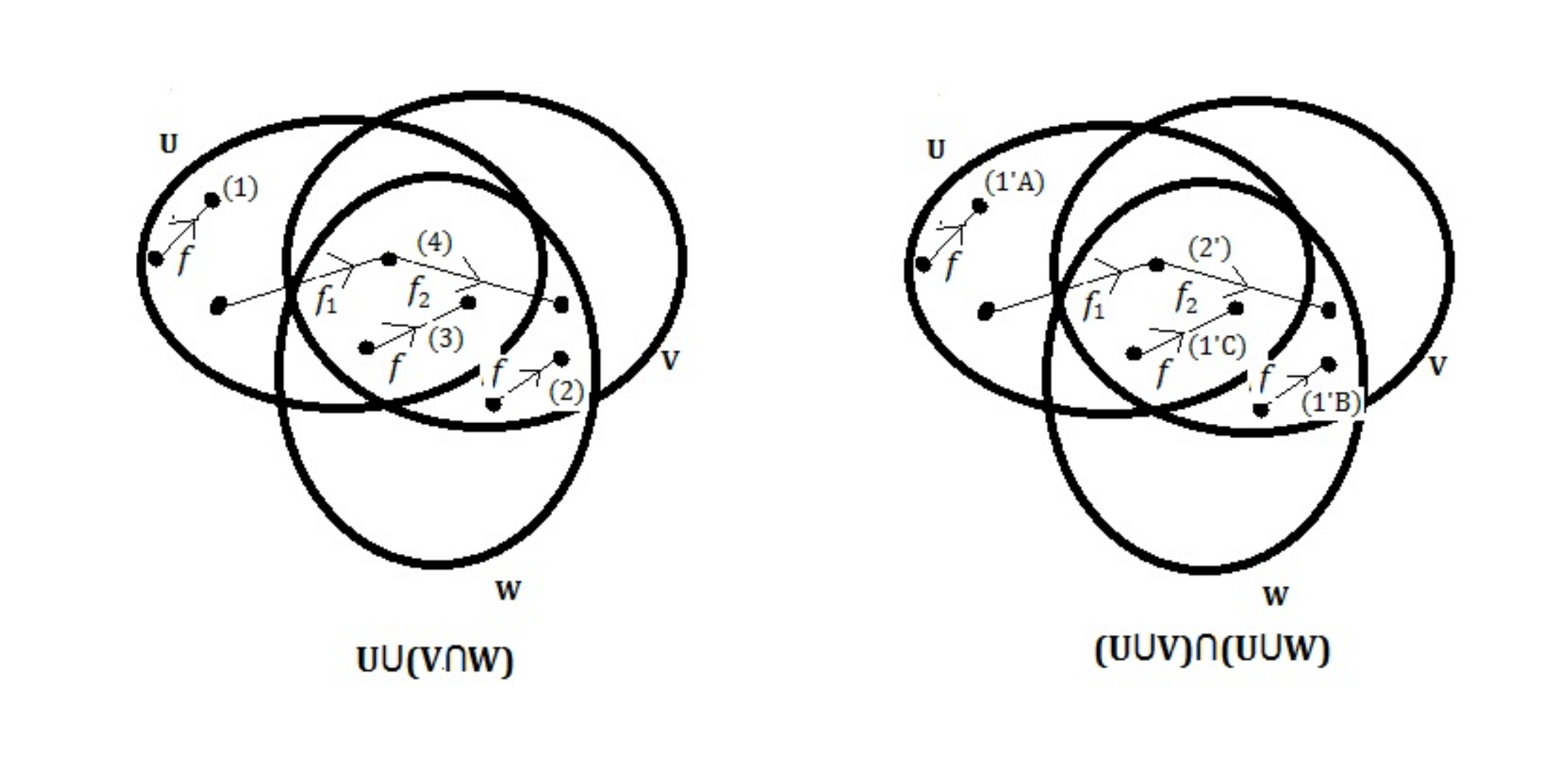}
\caption{ Comparison between morphisms in  $\mathbf U\cup({\mathbf V}\cap {\mathbf W})$ and $(\mathbf U\cup {\mathbf V})\cap (\mathbf U\cup{\mathbf W})$   }
\end{center}
\end{figure}\label{f:unionitii}

As before, we classify the morphisms in the left hand  and right hand sides of (ii) and compare. 

Let $f\in \Mor\biggl(\mathbf U\cup({\mathbf V}\cap {\mathbf W})\biggr).$ Then following are the possibilities (see the left hand side of Figure 5).
\begin{itemize}
\item[]{ (1) $f\in \Mor(\mathbf U).$
}
\item[]{ (2) $f\in \Mor({\mathbf V}\cap {\mathbf W})=\Mor({\mathbf V})\cap\Mor({\mathbf W}).$
}
\item[]{ (3) $f\in \Mor(\mathbf U)\cap\Mor({\mathbf V}\cap {\mathbf W})=\Mor({\mathbf U})\cap\Mor({\mathbf V})\cap\Mor({\mathbf W})=\Mor(\mathbf U\cap{\mathbf V}\cap {\mathbf W}).$
}
\item[]{ (4) $f=f_2\circ f_1,$ where $f_2, f_1\in \Mor(\mathbf U)\cup\Mor(\mathbf V\cap \mathbf W)$ and $s(f_2)=t(f_1)\in \Obj(\mathbf U)\cap \Obj(\mathbf V\cap \mathbf W)=\Obj(\mathbf U)\cap \Obj(\mathbf V)\cap \Obj(\mathbf W)=\Obj(\mathbf U\cap\mathbf V\cap \mathbf W).$
}
\end{itemize}
On the other hand, if \\
$f\in \Mor\biggl((\mathbf U\cup {\mathbf V})\cap (\mathbf U\cup{\mathbf W})\biggr)=\Mor\biggl(\mathbf U\cup {\mathbf V}\biggr)\bigcap \Mor\biggl(\mathbf U\cup{\mathbf W}\biggr),$ then possibilities are as follows (see the right hand side of Figure 5).
\begin{itemize}
\item[]{ (\'1) $f\in \biggl(\Mor(\mathbf U)\cup \Mor(\mathbf V)\biggr)\cap\biggl(\Mor(\mathbf U)\cup \Mor(\mathbf W)\biggr).$ Since 
\begin{equation}\nonumber
\begin{split}
&\Mor(\mathbf U)\cup \Mor(\mathbf V)\subset \Mor\biggl(\mathbf U\cup {\mathbf V}\biggr),\\
&\Mor(\mathbf U)\cup \Mor(\mathbf W)\subset \Mor\biggl(\mathbf U\cup {\mathbf W}\biggr)\hbox{[by \eqref{E:spanmorrewri}]},
\end{split}
\end{equation} 
we have
\begin{equation}\nonumber
\begin{split}
&\biggl(\Mor(\mathbf U)\cup \Mor(\mathbf V)\biggr)\bigcap \biggl(\Mor(\mathbf U)\cup \Mor(\mathbf W)\biggr)\subset  \Mor\biggl(\mathbf U\cup {\mathbf V}\biggr)\bigcap \Mor\biggl(\mathbf U\cup{\mathbf W}\biggr),\\
&\Rightarrow \Mor(\mathbf U)\cup\biggl(\Mor(\mathbf V)\cap \Mor(\mathbf W)\biggr)\subset \Mor\biggl(\mathbf U\cup {\mathbf V}\biggr)\bigcap \Mor\biggl(\mathbf U\cup{\mathbf W}\biggr)\\
\end{split}.
\end{equation} 
So we further classify as
\begin{itemize}
\item[]{(\'1A) $f\in \Mor(\mathbf U).$
}
\item[]{(\'1B) $f\in\biggl(\Mor(\mathbf V)\cap \Mor(\mathbf W)\biggr)=\Mor\biggl(\mathbf V\cap \mathbf W\biggr).$
}
\item[]{(\'1C) $f\in\Mor(\mathbf U)\cap\biggl(\Mor(\mathbf V)\cap \Mor(\mathbf W)\biggr)=\Mor\biggl(\mathbf U\cap \mathbf V\cap \mathbf W\biggr).$
}
\end{itemize}
}

\item[]{ (\'2) $f=f_2\circ f_1,$ where $f_2, f_1\in \Mor(\mathbf U) \cup\Mor(\mathbf V\cap \mathbf W)$ and $s(f_2)=t(f_1)\in \biggl(\Obj(\mathbf U)\cap \Obj(\mathbf V)\biggr)\bigcap \biggl(\Obj(\mathbf V)\cap \Obj(\mathbf W)\biggr)=\Obj(\mathbf U\cap\mathbf V\cap \mathbf W).$
}
\end{itemize} 

We make following correspondence between two sets of classifications:
\begin{equation}\nonumber
\begin{split}
&\rm{(1)}\Longleftrightarrow \hbox{(\'1A)}\\
&\rm{(2)}\Longleftrightarrow \hbox{(\'1B)}\\
&\rm{(3)}\Longleftrightarrow \hbox{(\'1C)}\\
&\rm{(4)}\Longleftrightarrow \hbox{(\'2)}
\end{split}
\end{equation}
and conclude
$${\mathbf U}\cup ({\mathbf V}\cap {\mathbf W})=({\mathbf U}\cup {\mathbf V})\cap({\mathbf U}\cup {\mathbf W}).$$
\subsection*{Proof of identity (iii)}
Same methodology  can be adopted to prove,
$${\mathbf U}\cap ({\mathbf V}\cup {\mathbf W})=({\mathbf U}\cap {\mathbf V})\cup({\mathbf U}\cap {\mathbf W}).\fbox{}$$
\section*{Concluding remarks}
In this paper we have developed a framework of  ${\mathbf {Cat}}$-valued sheaves over a category ${\widetilde {\mathcal O}}(\mathbf B)$ of subcategories of a  topological groupoid  ${\mathbf B}.$ Our starting point was the definition of ${\mathbf {Cat}}$-valued presheaf  introduced in {}[13]. We have constructed the ${\mathbf {Cat}}$-valued sheaf of local functorial sections on $\mathbf B$ for a fixed category $\mathbf C$. We have further shown that if we replace ${\mathbf C}$   by a categorical group $\mathcal G,$ we obtain a ${\mathbf {CatGrp}}$-valued sheaf. In traditional sheaf theory, sheaf of sections on a given topological space is the basic object of interest. In fact, every sheaf defined on a topological space can be realized as a sheaf of sections on the so called \textit{\'etal\'e space} corresponding to the given topological space [27]. It would be interesting to see if similar consideration also arises in the context of this paper; that is, whether there exists an ``{\'etal\'e category}", such that any ${\mathbf {Cat}}$-valued sheaf on $\mathbf B$ can be realized as a sheaf of functorial sections on the ``{\'etal\'e category}".

Lastly, we note that it would have been more natural if we had defined ${\mathbf {Cat}}$-valued presheaf to be a $2$-functor
$${\mathcal C}^{\rm op}\longrightarrow {\mathbf {Cat}}.$$
Because, ${\mathbf {Cat}}, {\mathcal C}, {\widetilde {\mathcal O}}(\mathbf B)$ are all $2$-categories. To reduce the complexity we have purposefully ignored the natural higher structures which fore mentioned categories posses. However, it is not very difficult task to extend our frame work to the higher level of enrichment.

\vskip .3 cm

 {\bf{ Acknowledgments.} } \textit{The author gratefully acknowledges   suggestions received  from Ambar N Sengupta and Amitabha Lahiri. The author would like to thank International Centre for Theoretical Sciences, TIFR, Bangalore for their kind hospitality, where part of this paper was written.}

\end{document}